\documentclass[paper,11pt]{AAS}	
\usepackage{xcolor}
\usepackage{bm}
\usepackage{amsmath}
\usepackage{subfigure}
\usepackage[colorlinks=true, pdfstartview=FitV, linkcolor=black, citecolor= black, urlcolor= black]{hyperref}
\usepackage{overcite}
\usepackage{footnpag}			      	
\usepackage{amsfonts}
\usepackage{array}
\usepackage{multirow}
\usepackage{adjustbox}

\usepackage{outlines}

\usepackage{grffile}

\DeclareGraphicsExtensions{.png,.pdf}
\graphicspath{{./figures/}}

\newcounter{ibid}
\newcounter{tempcnt}
\newcommand{\thanksibid}[1]{%
  \begingroup
  \setcounter{tempcnt}{\value{footnote}}%
  \setcounter{ibid}{#1}%
  \def\thefootnote{\fnsymbol{ibid}}%
  \footnotemark%
  \setcounter{footnote}{\value{tempcnt}}%
  \endgroup
}

\PaperNumber{}

\begin{document}

\title{
Isolating Neighborhood Trajectory Computations in Non-Autonomous Systems
Including the Elliptic Restricted Three-Body Problem
}

\author{Rodney L. Anderson,\thanks{Technologist,
Jet Propulsion Laboratory,
California Institute of Technology,
4800 Oak Grove Drive,
M/S 301-121,
Pasadena, CA 91109
}
\hspace{2pt}
 Robert W. Easton,\thanks{Professor Emeritus,
Applied Mathematics,
University of Colorado at Boulder,
Boulder, CO 80309
\newline
\copyright \hspace{2pt} 2023 California Institute of Technology. Government sponsorship acknowledged.
}
\hspace{2pt} and Martin W. Lo\thanksibid{1}
}

\maketitle{} 		

\begin{abstract}

Isolating block and isolating neighborhood methods have previously been
implemented to find transit trajectories and orbits around libration points in
the autonomous circular restricted three-body problem.  For some applications,
the direct computation of these types of trajectories in non-autonomous models
more closely approximating real-world ephemerides is beneficial.  Here, we
apply isolating neighborhood methods to non-autonomous systems, including the
elliptic restricted three-body problem (ERTBP).  Specifically, simplified
isolating neighborhood boundaries are computed around libration points in the
ERTBP.  These boundaries are used in combination with a bisection method to
compute the forward asymptotic trajectories of the isolated invariant set and
track orbits around a libration point.

\end{abstract}

\section{Lead Paragraph}

The computation of transit trajectories and orbits around libration points in
the autonomous circular restricted three-body problem (CRTBP) has previously
been enabled with the use of isolating blocks and isolating neighborhoods.  The
implementation of these methods in non-autonomous systems allows for a more
straight-forward application of these algorithms to real-world models, although
their use with non-autonomous systems presents several challenges.  We address
these issues by implementing new techniques in the non-autonomous elliptic
restricted three-body problem.  Specifically, the required boundaries are
tested across a range of energies that span the potential energies that
trajectories used in the computations may have, and the orbit tracking
algorithm is adapted to perform corrections at any required time.  The use of
these methods allows the tracking of quasiperiodic orbits having the
characteristics of 2-tori and 3-tori in the ERTBP.


\section{Introduction}

Isolating block and isolating neighborhood methods have previously been used to
describe the dynamics near libration points in the circular restricted
three-body problem (CRTBP).\cite{Conley:1965, Easton:1967, Conley:1968,
Easton:1979} Isolating blocks have been used to compute forward asymptotic
trajectories approaching the L$_2$ isolated invariant set and closely track
periodic and quasiperiodic orbits within the isolated invariant
set.\cite{Anderson:2017, Anderson:2020}  Isolating neighborhoods were then used
to simplify the computations and compute trajectories closely tracking
periodic and quasiperiodic orbits at higher energies including the halo orbit
and quasi-halo orbit families.\cite{Anderson:2019b}  Another approach was
introduced to compute desired periodic and quasiperiodic orbits within the
L$_2$ Lyapunov orbit in a surface of section directly without using asymptotic
approach trajectories.\cite{Anderson:2020b} This method also allows the
computation of trajectories within the chaotic region in this surface of
section.\cite{Anderson:2021e}  While the ability to compute these asymptotic
trajectories in the autonomous CRTBP is very useful,
our recent analyses have shown that it is possible to extend the isolating
neighborhood methods used for the CRTBP to non-autonomous systems.

In this study, we extend our methods, previously developed for use in the
autonomous CRTBP, to non-autonomous problems, including the elliptic restricted
three-body problem (ERTBP).
%
%
The isolating neighborhood boundaries developed for use in the CRTBP may be
used as a starting point for computing analogous isolating neighborhood
boundaries in the ERTBP.  In the CRTBP, a chosen boundary may be checked to
determine whether it is an isolating block boundary for a selected region by
computing tangent trajectories at a particular Jacobi constant and integrating
them forward and backward in time.\cite{Anderson:2019b}
In the ERTBP, there is no longer a Jacobi constant, but a similar boundary may
be computed around regions such as the L$_2$ point.  

In order to find an isolating neighborhood, we first specify a closed set $B$
in phase space around a region of interest, in this case the region around the
L$_2$ point. 
The closed set $B$ is a manifold with two boundary components - a ``right" and a
``left" boundary. We verify the exit behavior of  orbits to check that those that are tangent to the
right boundary exit through the right boundary and orbits tangent to the left
boundary exit left. 
%
%
%
In
this case, 
the isolating neighborhood boundary for the ERTBP must be checked at various epochs to ensure
that it remains an isolating neighborhood boundary for each epoch.
While checking the tangent trajectories at each epoch, velocities must be computed
for the tangent trajectories.  Since the Jacobi constant does not exist in the ERTBP,
a range of velocity magnitudes may be used to provide
%
evidence that this boundary still works as an isolating boundary
in the non-autonomous system.

Our approach using isolating neighborhoods serves as an alternative to standard
approaches that typically use some combination of expansions and differential
correction to compute quasiperiodic orbits in multibody systems.
In the autonomous CRTBP,
methods such as 
%
Fourier approaches,\cite{Gomez:2001f,Mondelo:2001} 
Lindstedt-Poincar\'e methods,\cite{Gomez:1998c} 
parameterization methods,\cite{Haro:2006,Haro:2007} 
normal form or semianalytical methods,\cite{Jorba:1998,Jorba:1999} 
Poincar\'e maps,\cite{Kolemen:2006,Kolemen:2012} and stroboscopic maps in
combination with
collocation\cite{Olikara:2010,Olikara:2010b,Olikara:2011,Olikara:2012} have
been used to compute quasiperiodic orbits.
Although much of the work in the non-autonomous ERTBP has focused
on the computation of periodic orbits,\cite{Broucke:1968, Broucke:1969, 
Bolotin:2005, Campagnola:2008, Cors:2001, Gin:2005,
Gomez:1986, Olle:1999, Koh:2016, Koh:2017c} some work has been done to explore
quasiperiodic tori in the ERTBP and other non-autonomous systems.
%
More specifically, Capinski and Zgliczy\'nski\cite{Capinski:2018} examined the
effect of perturbations in the planar ERTBP on Lyapunov orbits for small
eccentricities and showed that these orbits are perturbed to quasiperiodic
invariant tori.  Kumar, Anderson, and de la Llave explored resonant tori in the
planar ERTBP.\cite{Kumar:2020, Kumar:2022} Farres and Jorba\cite{Farres:2012}
explored orbits within the elliptic Hill problem, and Jorba, Jorba-Cus\'co, and
Rosales\cite{Jorba:2020, Rosales:2021} and Jorba and Villanueva
1997\cite{Jorba:1997b} studied the bicircular problem.  More recently,
F\'ernandez, Haro, and Mondelo studied quasiperiodic orbits in the
ERTBP.\cite{Fernandez:2022}

In this paper, the isolating neighborhood boundaries are computed for a range
of epochs and velocities.  The implementation of bisection methods for the
computation of the asymptotic trajectories approaching the isolated invariant
set in the ERTBP is also discussed. To provide insight into these trajectories,
we introduce several ERTBP coordinate systems.  We first describe an ERTBP
model using the standard constantly rotating coordinate system that is used to
study the CRTBP.\cite{Easton:2021}  We also use a non-uniformly rotating system
with an axis that remains aligned with the primaries,\cite{Hiday:1992b,
Hiday-Johnston:1994, Howell:1994} and finally, the commonly used ERTBP system
with a non-uniformly rotating and pulsating frame.\cite{Szebehely:1967}
Tracking methods for following periodic and quasiperiodic orbits within the
isolated invariant set are also developed and used to compute 
quasiperiodic orbits in the ERTBP.  The resulting orbits are described and
compared to results from the CRTBP.

%
%


%

%


%


%
%

\section{Models}

The previous analyses using isolating blocks and isolating neighborhoods have
focused on computations within the CRTBP.  Once the primaries are allowed to
move in orbits of non-zero eccentricity, the complexity of the problem
increases significantly, and the use of different ERTBP models is useful to aid
in both computations and visualization of the results.  Four different models
are introduced and described here.  The standard CRTBP is described first,
followed by an ERTBP model with a constantly rotating frame, an ERTBP model
with a non-uniformly rotating coordinate frame, and finally an ERTBP model with
both a non-uniformly rotating and pulsating coordinate frame.


\subsection{Circular Restricted Three-Body Problem}

While the primary focus of this work is on computing orbits within the ERTBP,
it is useful to compare the results to the orbits computed within the CRTBP.
In the CRTBP, the motion of a point mass is computed in a rotating frame
aligned with the primary and secondary masses which are constrained to move on
circular orbits.  The primary mass (the Earth for this study) is located on the rotating $x$ axis at
$E_{em} = \left( - \mu,0,0 \right)$, while the secondary mass (the Moon) is
located on the $x$ axis at $M_{em} = \left( \lambda,0,0 \right)$ where
$\lambda = 1 - \mu$.  The dimensionless mass $\mu = m_2/(m_1+m_2)$ is defined
where $m_1$ is the mass of the primary, and $m_2$ is the mass of the secondary.
The equations of motion in this model are 
%
\textcolor{black}{
\begin{equation}
\begin{split}
\ddot{x} &= \partial_{x}\Phi\left( x,y,z \right) \textcolor{black}{+} 2\dot{y} \\
\ddot{y} &= \partial_{y}\Phi\left( x,y,z \right) \textcolor{black}{-} 2\dot{x} \\
\ddot{z} &= \partial_{z}\Phi\left( x,y,z \right)
\end{split}
\end{equation}
}
where
\begin{equation}
\Phi\left( x,y,z \right) = \frac{1}{2}(x^{2} + y^{2}) + U(x,y,z)
\end{equation}
and
\begin{equation}
U\left( x,y,z \right) = \lambda/r_{1}\left( x,y,z \right) + \mu{/r}_{2}\left( x,y,z \right)
.
\end{equation}
Computing the partials gives
\begin{equation}
  \begin{aligned}
  \ddot{x} &=
2 \dot{y} 
+ 
x
-\frac{\left(x + \mu\right)\,\left(1 - \mu\right)}
{r_1^3}
-\frac{\mu \,\left(x -1 + \mu\right)}
{r_2^3}
\\
  \ddot{y} &= -2\dot{x} 
+ y
\left[
1
- \frac{1 - \mu}{r_1^3}
- \frac{\mu} {r_2^3}
\right] \\
  \ddot{z} &= 
z
 \left[ -\frac{1-\mu}{r_1^3} -\frac{\mu}{r_2^3} \right]
.
  \end{aligned}
\end{equation}
The distance between the point mass and the primary is $r_1 \left(x,y,z \right)$, and the distance
from the point mass to the secondary is $r_2 \left( x,y,z \right)$.
This study focused on the Earth-Moon system using a mass ratio
of $\mu = 1.2150584270571545 \times 10^{-2}$.
A constant of motion (the Jacobi constant) exists in the CRTBP, and it is defined as $C = -2J$
where
\begin{equation}
J = \frac{1}{2}\left\langle \dot{q},\dot{q} \right\rangle -  \Phi(q)
\end{equation}
or
\begin{equation}
\label{jacobiconstant}
C = -\dot{x}^2 - \dot{y}^2 - \dot{z}^2 + x^2 + y^2 + \frac{2(1-\mu)}{r_1} + \frac{2 \mu}{r_2}
.
\end{equation}

\subsection{Elliptic Restricted Three-Body Problem}


One of the most common formulations of the ERTBP\cite{Szebehely:1967} uses a
non-uniformly rotating coordinate system where the $x$ axis is aligned with the
primaries, and the length unit is chosen to be the varying distance between the
primary and the secondary.  
The foundation for deriving the equations of motion for this system is laid here
by first deriving the equations of motion in both uniformly and non-uniformly rotating coordinate systems.
A new derivation for the non-uniformly rotating and pulsating coordinate system
is then described.  This results in the same equations of motion obtained using Szebehely's\cite{Szebehely:1967}
derivation based on the work by Scheibner,\cite{Scheibner:1866} Petr and
Nechvile,\cite{Petr:1918} Nechvile,\cite{Nechvile:1926} and
Rein.\cite{Rein:1940}

In the following, capital letters are used for inertial coordinates. 
Starting with the full three-body problem and making the third mass zero gives
the equations of motion for a restricted three-body problem as
%
\begin{equation} \label{nonuniformeq3}
\ddot Q_1 = G m_2 |Q_2 -Q_1|^{-3} (Q_2 -Q_1 )
\end{equation}
\begin{equation} \label{nonuniformeq4}
\ddot Q_2 = G m_1 |Q_2 -Q_1|^{-3} (Q_1 -Q_2 )
\end{equation}
\begin{equation} \label{nonuniformeq5}
\ddot Q_3 = G m_1 |Q_1 -Q_3|^{-3} (Q_1 -Q_3 ) +
G m_2 |Q_2 -Q_3|^{-3} (Q_2 -Q_3 )
.
\end{equation}
In these equations, $Q$ is position, and $V$ is velocity.  The subscripts 1, 2,
and 3 correspond to the primary, the secondary, and the infinitesimal mass,
respectively.
G is the universal
gravitional constant, and the masses $m_1$ and $m_2$ are defined as they are in
the CRTBP.
Next, set the center of mass $ m_1 Q_1 + m_2 Q_2 = 0$.
Now, set $K=Gm_1+Gm_2$ and define $ \mu = G m_2/K$ as in the CRTBP.
Equations \ref{nonuniformeq3} and \ref{nonuniformeq4} decouple from Equation \ref{nonuniformeq5} and can be
combined as Kepler's equation:
\begin{equation} 
\label{keplereq4}
 \ddot Q^*= -K |Q^*|^{-3}Q^*
.
\end{equation}
If $Q^*(t)$ is a solution of Equation \ref{keplereq4} then
$Q^*_1(t)= -\mu Q^*(t)$ and $Q^*_2(t) =
(1-\mu) Q^*(t) $ are solutions of Equations \ref{nonuniformeq3} and
\ref{nonuniformeq4}.

To derive a model for the ERTBP one chooses the
eccentricity $e$ and the semi-major axis $a$ for the ellipse $ r(\nu) =
\frac{a(1-e^2)}{1+e \cos(\nu)}$.  Then one solves for initial conditions for a
solution $Q^*(t)$ of Equation \ref{keplereq4} that traverses the ellipse. Set $Q_1^*(t)= \mu Q^*(t)$
and $Q_2^*(t)= (1-\mu) Q^*(t)$. Then with $Q_3$ replaced by  $Q$, the elliptic
model is given by the equation 
\begin{equation} 
\label{keplereq5}
\ddot Q = G m_1 |Q_1^* -Q|^{-3} (Q_1^* -Q ) +
G m_2 |Q_2^* -Q|^{-3} (Q_2^* -Q )
.
\end{equation}
 
 Elliptic parameters for a solution of Kepler's Equation \ref{keplereq4} are
related to the semi-major axis $a$, the energy $E$, and angular momentum
$\sigma$ of the solution. Using the relations $e = \sqrt {1+2E \sigma^2/K^2}$
and $a(1-e^2)=\sigma^2/K$ one can choose initial conditions for a solution
$(Q^*,\dot Q^*)$
of Equation \ref{keplereq5} that traverses the ellipse. We choose the initial conditions
$Q^*(0) = [a(1-e);0,0]$ and 
 $\dot Q^*(0) = [0; \sqrt \frac{K(1+e)}{a(1-e)};0]$.
%
A normalized distance unit can be set so that $a(1-e) = 1$. For this unit, 
 $Q^*(0)=[1;0;0]$ and $ \dot Q^*(0)=[0;\sqrt{K(1+e)};0]$. Alternately one may set $a = 1$. 
%
%
 The CRTBP model results from the choice
 $e = 0 $. A hyperbolic ``fly-by" model results from the choice $ e > 1$. 
%

\subsection{Uniformly Rotating ERTBP Coordinate System}

A formulation of the ERTBP using a coordinate frame with fixed rotation was
introduced by Easton.\cite{Easton:2021} This coordinate frame allows an easier
comparison with the CRTBP since it uses the same frame. Both the primary and
secondary will move in position in this frame.  We will refer to Easton for the
detailed derivation, but an overview is given here.

The CRTBP is derived when a circular solution $Q^*(t)$ of radius $a$ of
Kepler's equation
is chosen, and a constantly rotating coordinate system is
introduced. This system has time, position, and velocity coordinates $(t,q,
\dot q)$ and the rotation matrix 
\begin{equation}
R(\omega t) =
\left(
\begin{array}{ccc}
  \cos (\omega t)& -\sin(\omega t)  & 0  \\
  \sin(\omega t)&  \cos (\omega t) &  0 \\
  0 & 0  &  1 
\end{array}
\right)
.
\end{equation}
The circular solution is $Q^* (t) = R(\omega t) q_0$, with 
$q_0 =(a,0,0)$ and 
$\omega = \sqrt{\mu / a^3 } $. 

The dynamics of the ERTBP in the uniformly rotating system require the primary
and secondary masses to move in position in both the $x$ and $y$ directions. 
%
In order to derive the equations of motion in this system, we first
introduce uniformly rotating coordinates $T_U: (t,q,\dot q) \rightarrow (t, Q, \dot Q) $.
Then we set $Q=R(\omega t)q$ and
\begin{equation}
\dot R = RA \omega
\end{equation} 
where
\begin{equation}
A=
\left(
\begin{array}{ccc}
 0& 1  & 0  \\
 -1&  0 &  0 \\
  0 & 0  &  0 \end{array}
\right)
.
\end{equation}
Now, 
\begin{equation} \label{e}
\ddot Q = R[\ddot q + 2\omega A \dot q + \omega^2 A^2 q]
.
\end{equation}
The right hand side of Equation \ref{keplereq5} is 
\begin{equation} \label{e}
 G m_1 |Q_1^* -Q|^{-3} (Q_1^* -Q ) +
G m_2 |Q_2^* -Q|^{-3} (Q_2^* -Q )
.
\end{equation}
The elliptic solution $Q_e^*(t)$ of Kepler's equation that we use is a perturbation of the circular solution $Q^*$ used for the CRTBP model.
\begin{equation}
Q_e^*(t) = R(\nu (t)) r(\nu(t)) Q^*_0 =
 R(\omega t)R(\nu (t)-\omega t) q_e^*(t)
\end{equation}
with $q_e^*(t) = r(\nu(t))Q^*_0$. 
\begin{equation}
Q_j^*- Q = R(\omega t)[R(\nu (t)-\omega t) (q_j^* -q)]
\end{equation} 
with $q_1^* = -\mu q_e^*$ and $q_2^* = (1-\mu) q_e^*$.
The rotation matrix $R(\nu (t)-\omega t)$ compensates for the difference between the true anomaly as a function of time and the constant rotation rate.
The elliptic model in uniformly rotating coordinates is 
\begin{equation} \label{f}
\ddot q + 2\omega A q + \omega^2 A^2 q = 
R(\nu (t)-\omega t)[G m_1 d_1^{-3} ( q_1^* -q) +
G m_2 d_2^{-3} ( q_2^* -q )]
.
\end{equation} 
The distances are
$ d_1 = |q_1^* -q|, d_2 = |q_2^* -q|$.

\subsection{Non-Uniformly Rotating ERTBP Coordinate System}
We will use the true anomaly $\nu(t)$ of the solution $Q^*(t)$ to non-uniformly
rotate coordinates and to keep the primary and secondary masses on the x-axis.
The true anomaly will replace time as the independent variable. We use
prime to denote differentiation with respect to the true
anomaly.
Note that even though time is replaced with true anomaly, the phase portrait of the system
is still preserved.
Refer to Appendix A for a description of an alternative formulation of the 
ERTBP in non-uniformly rotating coordinates using time as the independent variable.

The true anomaly of the solution $Q^*(t)$ satisfies the equation 
\begin{equation} \label{eq9}
  \dot \nu = \sigma / r^2(\nu)
.
\end{equation}
First, we define a rotation matrix
\begin{equation}
R(\nu) =
\left(
\begin{array}{ccc}
  \cos (\nu)& -\sin(\nu)  & 0  \\
  \sin(\nu)&  \cos (\nu) &  0 \\
  0 & 0  &  1 
\end{array}
\right)
.
\end{equation}
Then $R^\prime(\nu)= R(\nu)A$ with
\begin{equation}
A=
\left(
\begin{array}{ccc}
 0& 1  & 0  \\
 -1&  0 &  0 \\
  0 & 0  &  0 \end{array}
\right)
.
\end{equation}
Next, introduce non-uniformly rotating coordinates $T_N: (t,X,Y) \rightarrow (t, Q, V) $, and
set $Q=R(\nu(t))X$.  Then 
$\dot R = R^\prime \dot \nu = R A \dot \nu , \ 
 \ 
\dot Q = R^\prime \dot \nu X + R \dot X
$. 
The true anomaly satisfies the conditions
\begin{equation} \label{eq10}
\dot \nu = \sigma r^{-2}, \nu(0)=0, \ 
\sigma = \sqrt {(Gm_1 + Gm_2)a(1-e^2)}
.
\end{equation}
The elliptic solution $Q^*(t)$ of Kepler's equation satisfies the equation
\begin{equation} \label{eq11}
Q^*(t) = R(\nu(t))r(\nu(t)) Q^*(0)
\end{equation}
where
\begin{equation} \label{eq12}
r(\nu) = \frac{a(1-e^2)}{\gamma(\nu)}, \ 
\gamma(\nu) = 1 + e \cos(\nu)
.
\end{equation}
In the non-uniformly rotating system the varying positions of the primary and secondary masses on the $x$-axis are 
$X_1^*(\nu(t))=-\mu  r(\nu(t)) Q^*(0)$ and
 $X_2^*(\nu(t))=(1-\mu) r(\nu(t))Q^*(0)$.
The equations for the ERTBP in non-uniformly rotating coordinates  $( \nu,X,Y)$ are derived as follows.
\begin{equation} \label{eq13}
\ddot Q = R[\ddot X + 2A \dot \nu \dot X
+A  \ddot \nu X + A^2 \dot \nu^2 X]
\end{equation}
 with 
$\dot X = X^\prime \dot \nu $ and 
$\ddot X = X^{\prime \prime} \dot \nu^2 + X^{\prime} \ddot \nu$.
This gives the result 
\begin{equation} \label{eq14}
\ddot Q = R\dot \nu^2 [X^{\prime \prime} +
2A X^\prime +
A^2 X ]+
R \ddot \nu [A X +  X^{\prime}] 
.
\end{equation}
Using the definition of $\nu$ we see that 
$\ddot \nu = (2 \gamma^\prime / \gamma) \dot \nu^2$ 
and $ \dot \nu^2 = \sigma^2r^{-4}$.
 Equation \ref{eq13} is now
 \begin{equation} \label{eq15}
 \ddot Q =R \sigma^2 r^{-4}[X^{\prime \prime} +
2A X^\prime +
A^2 X +
(2 \gamma^\prime / \gamma) (A X 
+  X^{\prime})]
\end{equation}
Express Equation \ref{e}  in terms of
$ (\nu,X)$: 
\begin{equation} \label{eq16}
F(\nu, X) = G m_1 \frac{ RX_1^*(\nu) -RX }{r_1^{3}}
+ G m_2  \frac{ RX_2^*(\nu) -RX}{r_2^{3}}
.
\end{equation}
We have $r_1 = |X_1^* - X| , r_2 = |X_2^* - X|$, 
$r \sigma^{-2} Gm_1  = \gamma^{-1} (1-\mu) , \ 
r \sigma^{-2} Gm_2  = \gamma^{-1}\mu$.
These expressions are used to give the result that
\begin{equation} \label{eq17}
F(\nu,X) =  R  \sigma^2 r^{-4}
\left[r^3\gamma^{-1}(1-\mu)\frac{ X_1^*(\nu) -X }{r_1^{3}}
+ r^3 \gamma^{-1} \mu \frac{ X_2^*(\nu) -X}{r_2^{3}} \right]
.
\end{equation}

The final result is the equation
\begin{equation} \label{eq18} 
 X^{\prime \prime} +
2A X^\prime +A^2 X +
(2 \gamma^\prime \gamma^{-1}) (A X +  X^{\prime}) = 
r^3 \gamma^{-1}(1-\mu) \frac{ X_1^*(\nu) -X }{r_1^{3}}
+r^3\gamma^{-1}\mu \frac{ X_2^*(\nu) -X}{r_2^{3}}
\end{equation}
with 
$ X_1^*(\nu) =  -\mu r(\nu)Q_0^*(0), 
\ 
 X_2^*(\nu) =  (1- \mu)r(\nu) Q_0^*(0)$.
 
 For the circular problem the term $ \gamma^{-1} \gamma^\prime$ is zero, and $r
= \gamma =1$. The true anomaly and the time variables are equal in this case.
 
\subsection{Non-Uniformly Rotating and Pulsating ERTBP Coordinate System }

The derivation of the equations of motion in non-uniformly rotating, pulsating
coordinates is new and uses the results for the non-uniformly rotating
coordinate system.
 The formulas that we use in the derivation are these: 
$ \dot \nu = \sigma r^{-2}$,
$\ddot \nu = 2 \phi \dot \nu^2$,
$\gamma = 1+ e \cos{(\nu)}$,
$\phi = \gamma^{-1} \gamma^\prime$,
$\phi^\prime = -\phi^2 + \gamma^{-1} \gamma^{\prime\prime}$, 
$r = p \gamma^{-1}$, 
$r^\prime = -r \phi$,
$r^{\prime\prime} = r(\phi^2 - \phi^\prime)$,
$R^\prime = RA$,
$\dot R = R^\prime \dot \nu$,
$\ddot R = RA \dot \nu^2 + RA \ddot  \nu$

Using the formulas we have 
$X = ru$,
$X^\prime = r(u^\prime - \phi u) $,
$X^{\prime\prime} =
 -r \phi (u^\prime - \phi u )+
  r(u^{\prime\prime} - \phi^\prime u - \phi u^\prime) $.
  Expanding 
  $X^{\prime\prime} =
 -r [-\phi u^\prime +  \phi^2 u +
  u^{\prime\prime} +
  \phi^2 u -\gamma^{-1} \gamma^{\prime\prime} u
 - \phi u^\prime] $.
After collecting terms we have  
  \begin{equation} \label{eq19}
  X^{\prime\prime} = r[u^{\prime\prime}+2\phi^2 u
  -2 \phi u^\prime 
  - \gamma^{-1} \gamma^{\prime\prime } u
 ] 
.
\end{equation}

A short but key calculation shows that
 $[1- \gamma^{-1}] = \gamma^{-1} \gamma^{\prime\prime}$.
Replacing each of the terms on the lefthand side of Equation \ref{eq18} and simplifying, the result is 
\begin{equation} \label{eq20} 
[X^{\prime \prime} +
2A X^\prime +A^2 X +
(2 \gamma^\prime \gamma^{-1}) (A X +  X^{\prime}) ] = 
r [u^{\prime\prime}+2A u^\prime+A^2u +(1-\gamma^{-1})u]
.
\end{equation}

\noindent Express the righthand side $F(\nu,X)$ of equation \ref{eq18}  in terms of
$ (\nu, u)$ as
\begin{equation} \label{eq21}
F(\nu, X) = r^3 \gamma^{-1}(1-\mu) \frac{ X_1^*(\nu) -X }{r_1^{3}}
+r^3\gamma^{-1}\mu \frac{ X_2^*(\nu) -X}{r_2^{3}}
\end{equation}
\begin{equation} \label{eq22}
F(\nu, u) =  r^4 \gamma^{-1}(1-\mu) \frac{ u_1^* -u }{r_1^{3}}
+r^4\gamma^{-1}\mu \frac{ u_2^* -u}{r_2^{3}}
\end{equation} 
with 
$ u_1^* =- \mu Q^*(0), \ u_2^* = (1-\mu)Q^*(0)$.
Note that $r_1 = r d_1 = r |u_1^*-u|$
 and update.
Then
\begin{equation} \label{eq23}
 F(\nu, u) = r [\gamma^{-1}(1-\mu) \frac{ u_1^* -u }{d_1^{3}}
+   \gamma^{-1}\mu \frac{ u_2^* -u}{d_2^{3}} ]
.
\end{equation}
The final result is 
\begin{equation} \label{eq24}
u^{\prime\prime}
+ 2   A  u^\prime
+   A^2 u +u
 =  (1 + e \cos(\nu))^{-1}[ u+(1- \mu) |u^*_1 -u|^{-3}
 (u^*_1 -  u)
+ ( \mu)|u^*_2 -u|^{-3}(u^*_2 - u) ]
.
\end{equation}
The result is formulated in vector-matrix notation as follows:
%
\begin{multline}
\begin{bmatrix}
x'' \\
y'' \\
z'' \\
\end{bmatrix}
+ 2
\begin{bmatrix}
0 & -1 & 0 \\
1 & 0 &0 \\
0 & 0 & 0  \\
\end{bmatrix}
\begin{bmatrix}
x' \\
y' \\
z' \\
\end{bmatrix}
+ \left(
\begin{bmatrix}
-1 & 0 & 0 \\
0 & -1 &0 \\
0 & 0 & 0  \\
\end{bmatrix}
\begin{bmatrix}
x \\
y \\
z \\
\end{bmatrix}
+
\begin{bmatrix}
x \\
y \\
z \\
\end{bmatrix}
\right)
= \\
\frac{1}{1 + e\cos(\nu)} \left[
\begin{bmatrix}
x \\
y \\
z \\
\end{bmatrix}
+
\frac{1 - \mu}{r_1^3}
\left(
\begin{bmatrix}
-\mu \\
0 \\
0 \\
\end{bmatrix}
 - 
\begin{bmatrix}
x \\
y \\
z \\
\end{bmatrix}
\right) -
\frac{\mu}{r_2^3}
\left(
\begin{bmatrix}
x \\
y \\
z \\
\end{bmatrix}
- 
\begin{bmatrix}
1-\mu \\
0 \\
0 \\
\end{bmatrix}
\right)
\right]
.
\end{multline}
%
The equations of motion for the system may then be written in the usual form as
\begin{equation}
\label{eomx}
\begin{aligned}
x'' 
&= 2y' + \frac{1}{1 + e\cos(\nu)} \left[ x
-\frac{1 - \mu}{r_1^3} (x + \mu)
- 
\frac{\mu}{r_2^3}
(x - 1 + \mu)
\right]
\end{aligned}
\end{equation}
\begin{equation}
\begin{aligned}
y'' 
&=
-2x' + \frac{y}{1 + e\cos(\nu)} \left[ 1 
-\frac{1 - \mu}{r_1^3}
-
\frac{\mu}{r_2^3}
\right]
\end{aligned}
\end{equation}
\begin{equation}
\label{eomz}
z'' = -z + \frac{z}{1 + e\cos(\nu)} \left[ 1
-\frac{1 - \mu}{r_1^3}
- 
\frac{\mu}{r_2^3}
\right]
.
\end{equation}
The primary and secondary are fixed in this coordinate frame at $(-\mu,0,0)$ and $(1-\mu,0,0)$, respectively.

%


%
%
%
%

\section{Coordinate Systems Summary}

\noindent It is often convenient to transform states and trajectories between the different ERTBP coordinate systems
that we have defined so far.  We will name the various coordinates as

Inertial Coordinates: Variables $(t, Q, \dot Q) \in R^1 \times R^3 \times R^3 $

Constant Rotating Coordinates: Variables $(t, X, \dot X) \in R^1 \times R^3 \times R^3 $

Variable Rotating Coordinates: Variables $(\nu, X, X^\prime) \in R^1 \times R^3 \times R^3 $

Rotating Pulsating Coordinates: Variables $(\nu, u, u^\prime) \in R^1 \times R^3 \times R^3 $

\noindent Now, we summarize the transformations between the different coordinate systems.

1. Constant Rotating to Inertial:
 $F: (t, X,\dot X) \rightarrow (t,Q, \dot Q) $

\begin{equation}
Q = R(\omega t)X, \ \dot Q = \dot R(\omega t)X + R(\omega t) \dot X
\end{equation}

\begin{equation}
R(t) = \begin{bmatrix}
\cos(\omega  t) & \sin(\omega  t) & 0 \\ 
-\sin(\omega  t) & \cos(\omega  t) &0 \\
0&0&1
\end{bmatrix}
\end{equation}

2. Variable Rotating to Constant Rotating:
 $F^*: (\nu, X, X^\prime) \rightarrow (t,X, X^\prime) $

Time is expressed as a function of the true anomaly by solving the differential equation 
\begin{equation}
\label{timenu}
\frac{dt }{ d\nu} = r^2(\nu)/ \sigma 
\end{equation}
with initial condition $ t(0) = 0$, 
$\sigma = \sqrt {(Gm_1 + Gm_2)a(1-e^2)}$,
 and $ r(\nu) = \frac{a(1-e^2)}{ 1 + e \cos(\nu)} $.

3. Pulsating to Variable Rotating: $G: (\nu,u,u^\prime) \rightarrow (\nu,X, X^\prime) $

\begin{equation}
X = u/r = \frac{1 + e \cos \nu }{1-e^2} u
\end{equation}
\begin{equation}
X^\prime = u^\prime / r - (u/r^2)  r^\prime 
\end{equation}


Finally, it is often convenient to convert between time and true anomaly for different
formulations of the ERTBP.  These conversions may be performed using 
Equation \ref{timenu}
or standard procedures which are given in Appendix B for convenience.

%
%
%
%
%

\section{Illustrative Example Using Event Functions, Polyhedra, and Bisection}

The theory of normal forms for Hamiltonian systems spans more than a century
with contributions from many great mathematicians. The analysis is beautiful,
complex, and difficult.  It requires years of research at the highest level to
understand this material. A recent publication by Jorba et al.\
applies state of the art normal form techniques to calculate quasi-periodic
solutions of the non-autonomous bicircular problem.\cite{Jorba:2020}
While these methods yield powerful results, the methods introduced here provide
an alternative that we will use to obtain similar results that apply to
spacecraft mission design.  These methods also lay the foundation
for analyses in more complicated non-autonomous systems such as the full
ephemeris model.  This model may be thought of as a small perturbation of the 
ERTBP in systems such as the Earth-Moon, Jupiter-Europa, or Saturn-Enceladus systems
 among others.

To investigate the ERTBP in the vicinity of the collinear Lagrange points we
use the topological concept of connectedness in combination with numerical
solutions of equations of motion and event functions. In this work computer
programs were used to compute exit times from regions of interest and to locate
solutions that do not exit. A theoretical foundation motivating the method is
based on the study of isolating blocks and the Conley index.\cite{Conley:1971} 


To illustrate our approach, we use a bisection method to find solutions of a
non-autonomous system of differential equations that start in a specified
region of phase space at a specified time and do not exit within a specified
long time interval. 
The type of equation we use is a time-dependent vector field on $R^n$:
\begin{equation} \label{eq1}
 \dot{x}=F(t,x), x(t_0)=x_0
.
 \end{equation}
 An event function $E_k$ is defined by using a base point $b_k$ and a unit normal vector $u_k$. 
\begin{equation} \label{eq2}
E_k(x; b_k, u_k)= (x-b_k)^T u_k
\end{equation}
The half space defined using the event function $E_k$ is the space $H_k = \{x: E_k(x; b_k, u_k) \ge 0 \}$. 
 A polyhedron $B$  can be viewed as the
intersection of half spaces. A boundary face of $B$ is the set $ \partial_k B = B \cap \{E_k = 0 \}$. 
 A convenient event function that signals the exit of a solution from the polyhedron $B$ is the minimum of the event functions defining the polyhedron. 
Our problem is to find an initial condition $(t_0 ,x_0)$  in $B$ so that the
solution to the initial value problem does not exit within the time interval
$[t_0, t_1]$. The behavior of the vector field in Equation \ref{eq1} on the boundary of $B$ may be sufficient to guarantee a solution. 
The specific differential equation that we use as an example is given by
\begin{equation} \label{eq3}
 \begin{bmatrix}
x_1^{'} \\
x_2^{'}
\end{bmatrix} =
\begin{bmatrix}
x_1   \\
-x_2
\end{bmatrix}
+ \epsilon \begin{bmatrix}
\cos(t)x_1+\sin(t)x_2  \\
-\sin(t)x_1+\cos(t)x_2 \\
\end{bmatrix}
.
\end{equation}
The problem is to find orbits that remain in the square 
$B = [-1,1] \times [-1, 1]$ for a long time. The square is defined
 as the intersection of four half-spaces given by four event functions.
The event functions are $E_R (x)= 1-x_1$, $E_L(x) =  -1+x_1 $, 
$E_U (x)= 1-x_2$, $E_D (x)= -1+x_2$. When one of these functions is zero this indicates a right, left, up, or down exit.

\begin{figure}[ht!]
\centering
\subfigure[t = 0] { \label{vec000}
\includegraphics[width=0.23\textwidth]{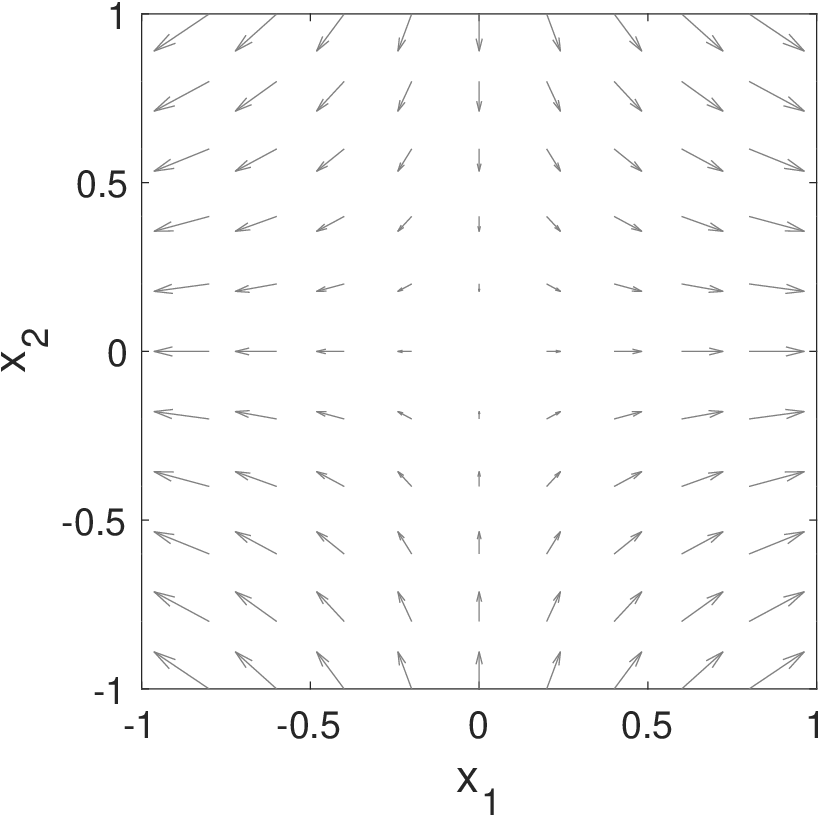}
} \subfigure[t = $\pi/2$] {\label{vec090}
\includegraphics[width=0.23\textwidth]{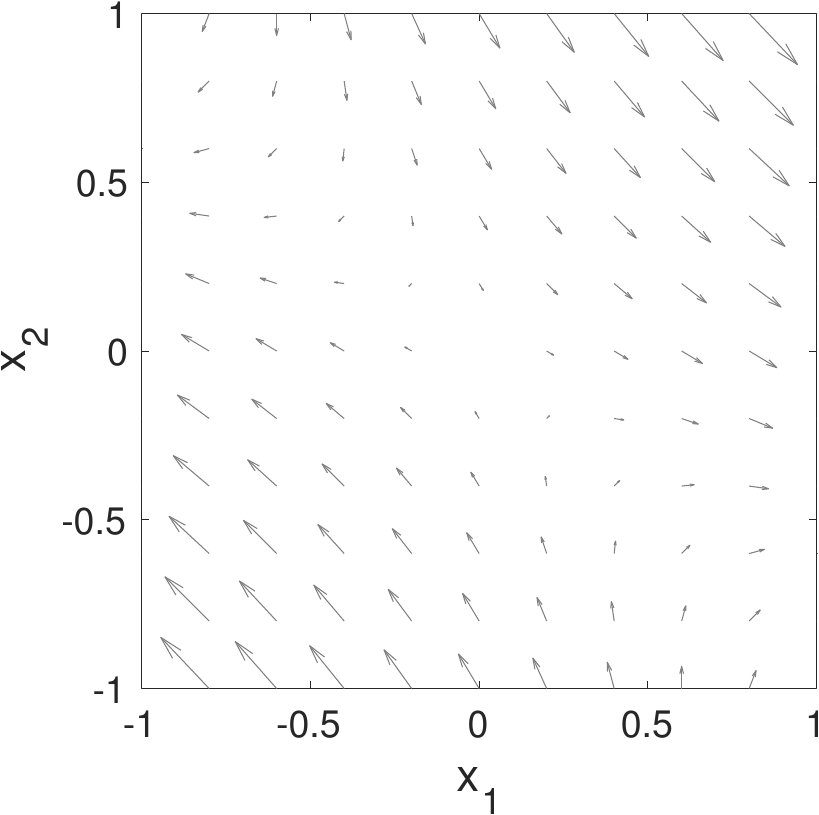}
}
\subfigure[t = $\pi$] {\label{vec180}
\includegraphics[width=0.23\textwidth]{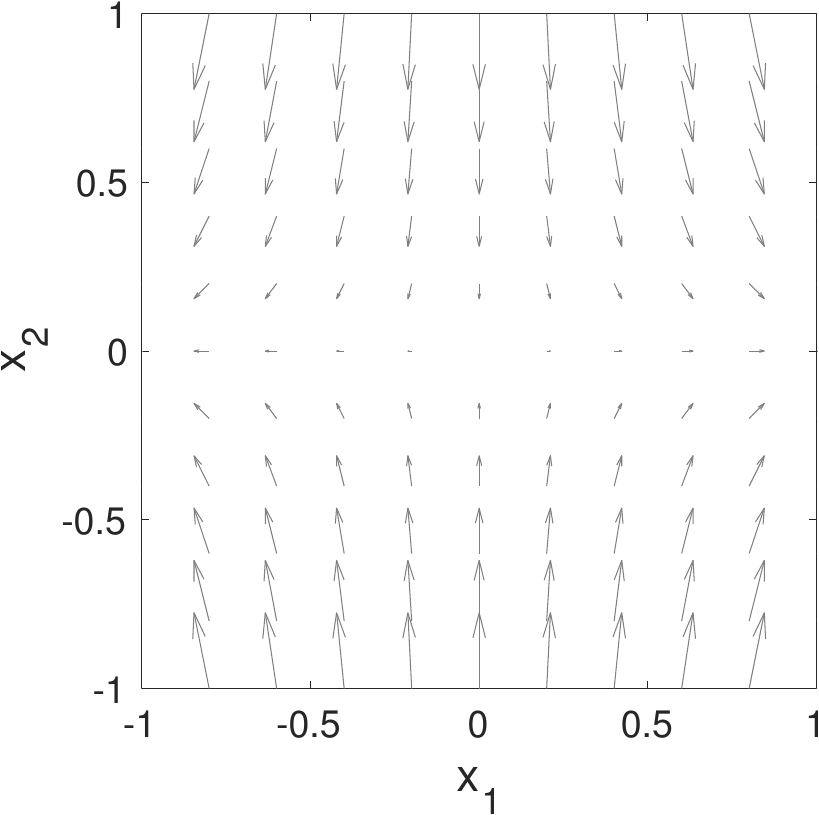}
}
\subfigure[t = $3\pi/2$] {\label{vec270}
\includegraphics[width=0.23\textwidth]{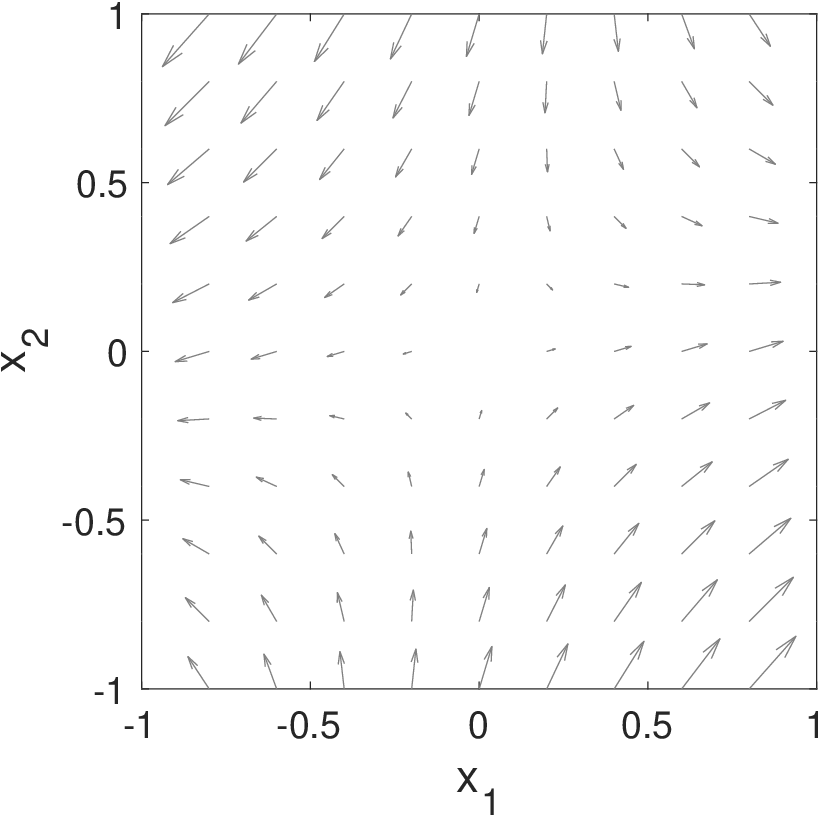}
}
\caption{
Vector fields for Equation \ref{eq3} computed for different times.
}
\label{vecfields}
\end{figure}

The vector field in Equation \ref{eq3} for $|\epsilon| < 1/\sqrt{2}$ is never tangent to
the boundary of $B$ and points outward on the right and left boundaries  and
inward at the top and bottom. This can be verified from the definition of the
vector field. The square $B$ is an isolating block with corners. In later, much
more complicated examples the exit behaviors of solutions tangent to boundaries
will be examined numerically. For this example the exit behavior of solutions
to the vector field with $ \epsilon = 0.6$ is illustrated for several initial
start times between 0 and $2\pi$ as shown in Figure \ref{vecfields}.

An orbit starting in $B$ and exiting must exit either through the right or left
boundary faces.  One can show that the set of initial conditions whose orbits
exit right form an open set and so do the initial conditions whose orbits exit
left. Given an initial time  and a line segment in $B$ with one endpoint
exiting right and the other exiting left, a topological argument shows that
there is a point $(x_1^*,x_2^*)$ on that segment such that the solution of
Equation \ref{eq3} remains in $B$ forward in time.  The topological reason is
that the line segment is a connected set, and it cannot be the union of
disjoint non-empty open sets. An approximation to the point $(x_1^*,x_2^*)$ can
be found by using a bisection method together with event functions and a
numerical method for solving the differential equation.  
It is worth noting that the bisection method will work for all small perturbations
of the basic linear equation in this example.  Normal form and parameterization
methods generally require more work to apply than the bisection approach.

\begin{figure}[ht!]
\centering
\subfigure[$t_0 = 0$] {\label{samplezero}
\includegraphics[width=0.48\textwidth]{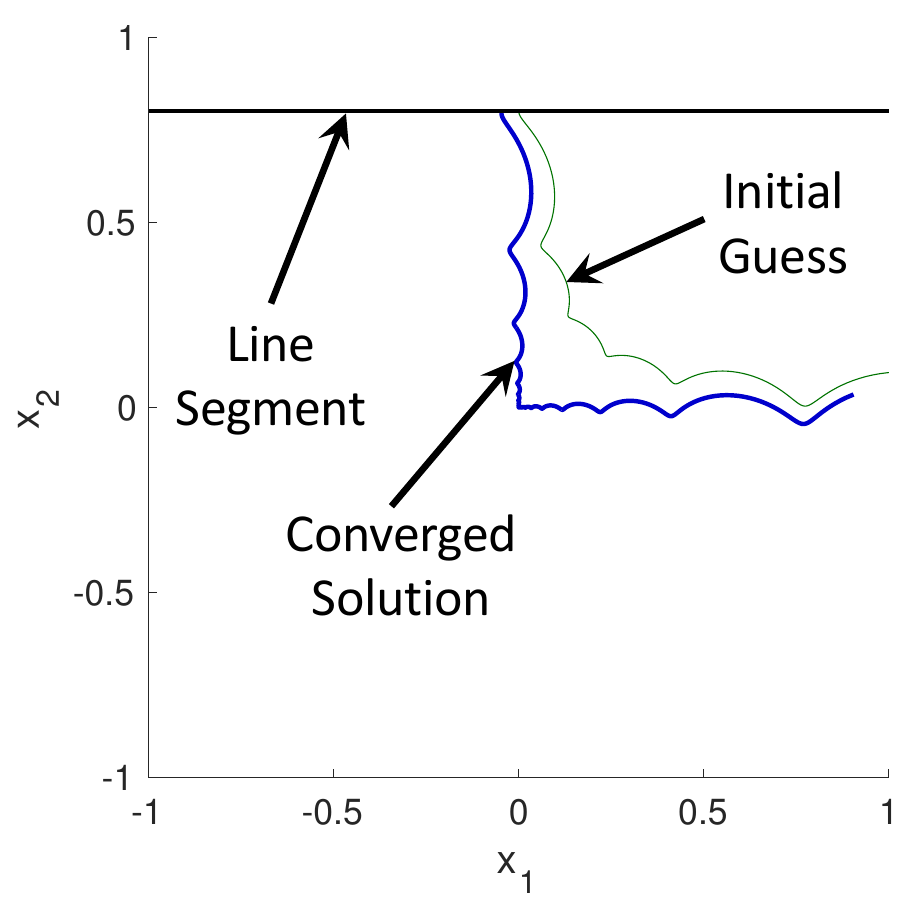}
}
\subfigure[$t_0 = \pi/2$] {\label{samplehalfpi}
\includegraphics[width=0.47\textwidth]{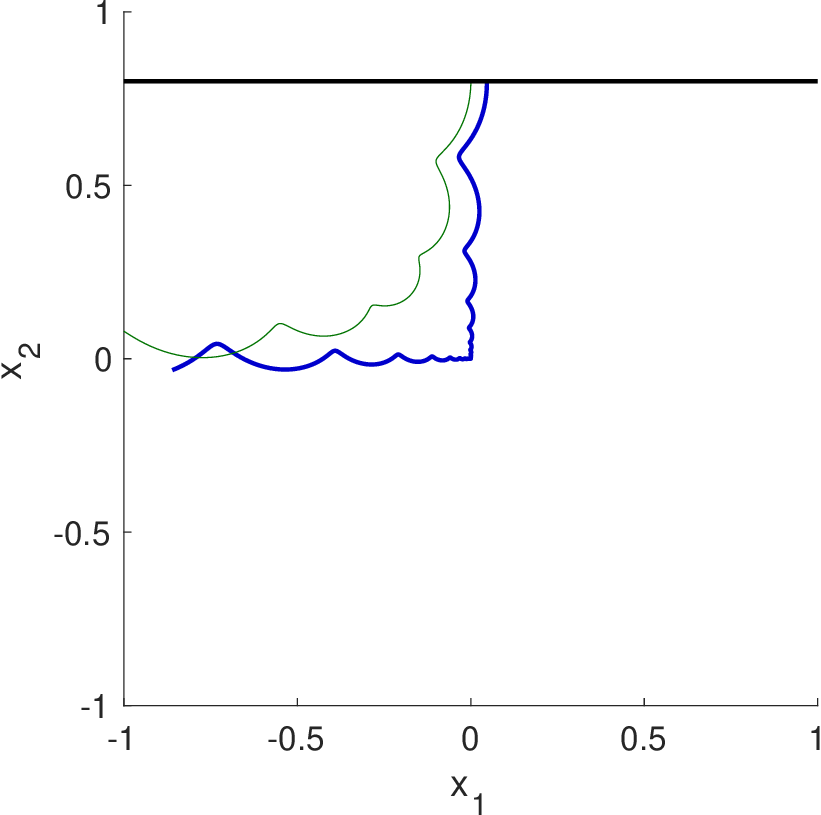}
}
\caption{
Converged solutions for $x_2 = 0.8$ at $t_0 = 0$ and $t_0 = \pi/2$.
}
\label{sampleboth}
\end{figure}

This problem was solved numerically, and the results for several cases are
described next.  As a first step, a line segment must be selected with end
points that exit left and right.  A simple case is chosen with $x_2 = 0.8$, and
$-1 \le x_1 \le 1$.  A quick check shows that the point $(-1,0.8)$ exits left, and $(1,0.8)$
exits right.  The next step is to select an initial time for the non-autonomous
system, and then a bisection method may be performed to find the value of $x_1$
that remains in the unit square as long as possible.  (Note that this
calculation is practically limited by the numerical precision of the computer.)
The results for initial times $t_0 = 0$ and $t_0 + \pi/2$ are shown in
Figure \ref{sampleboth}.
An initial guess at $(0,0.8)$ is shown in green in each plot, and the converged
final solution from the bisection method is shown in blue.  The converged solutions
in these cases stays within the unit square for a very long time in each of these
cases, but because of the limits of numerical precision of the computer, they will
eventually start to wander away.  The converged solutions do stay within the unit square
much longer than the initial guess from the bisection method though.

\section{Isolating Block/Neighborhood Boundaries in the CRTBP and ERTBP}

In previous papers, we have explored the use of both isolating block and
isolating neighborhood boundaries in the CRTBP for a variety of
systems.\cite{Anderson:2017b,Anderson:2019b,Anderson:2020,Anderson:2020b,Anderson:2021e}
Isolating block boundaries require that tangent trajectories immediately leave,
or ``bounce off" the boundaries.  Isolating neighborhood boundaries are less
stringent in the requirement for the test on these tangent trajectories.  In
this case, given ``left" and ``right" exit boundaries obtained in combination
with the Hill's region, the requirement is that tangent trajectories on the
left boundary exit left, and tangent trajectories on the right boundary exit
right.  We focus on the use of isolating neighborhood boundaries here.



\begin{figure}[ht!]
\centering
\subfigure[Left Exit] {\label{crtbpinnerin}
\includegraphics[width=0.48\textwidth]{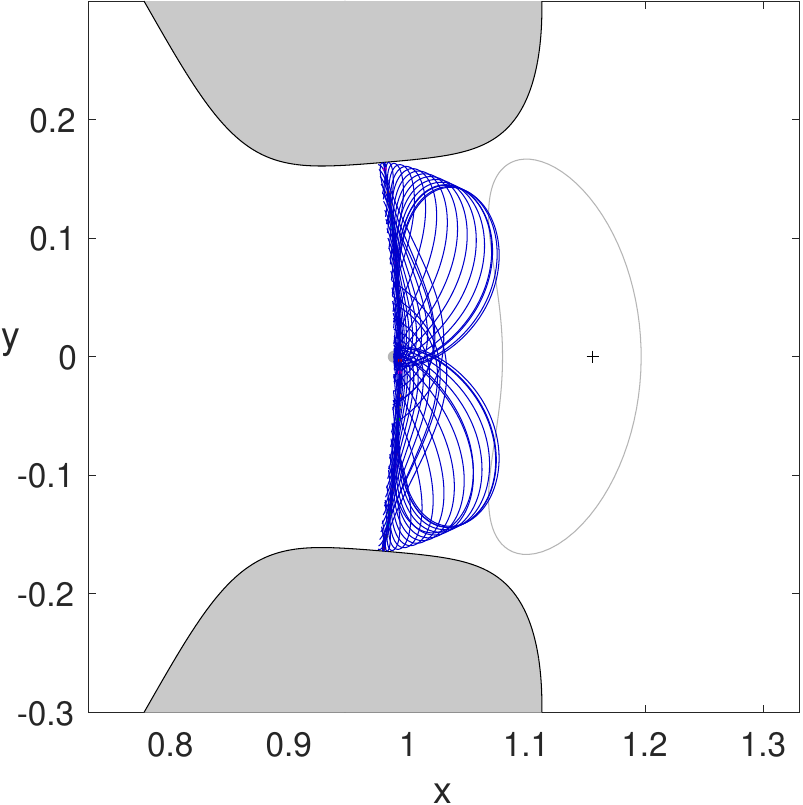}
}
\subfigure[Right Exit] {\label{crtbpouterin}
\includegraphics[width=0.48\textwidth]{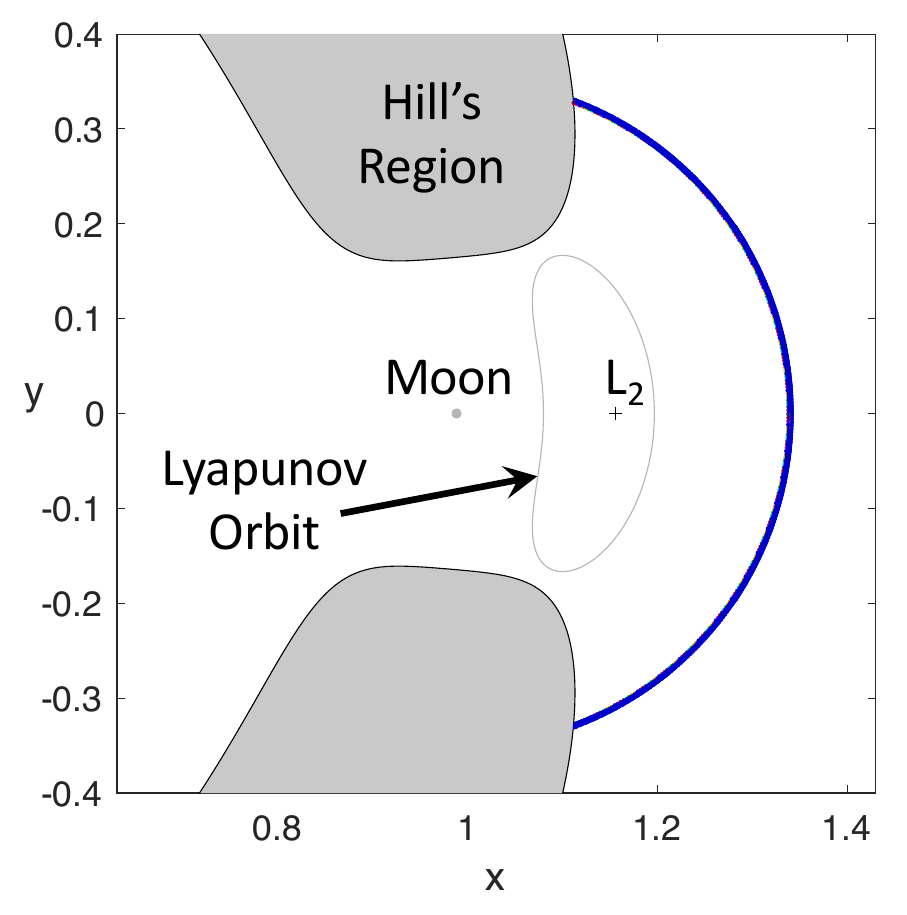}
}
\caption{
Test showing the tangent trajectories exiting the isolating neighborhood in the CRTBP.
The L$_2$ Lyapunov orbit at C = 3.10 is shown in gray for comparison with the isolating
neighborhood boundaries.
}
\label{crtbpneighborhoodtest}
\end{figure}

Cylinders were used as isolating neighborhood boundaries in Anderson, Easton,
and Lo,\cite{Anderson:2020b}  and a similar approach is used here.  As a first
step, isolating neighborhood boundaries are defined and checked for the CRTBP.
In this case, the inner boundary $\partial_L N$ is defined by $r_L^2 = x^2 +
y^2$ where $r_L = 0.993$, and the outer boundary $\partial_R N$ is defined by
$r_R^2 = (x -1 + \mu)^2 + y^2$ where $r_R = 0.35$.  
%
%
%
In each case, the validity of the
boundaries is verified for the planar problem by integrating tangent
trajectories both forward and backward until they exit the isolating
neighborhood.  
If tangent trajectories on the grid on the right side
exit right and tangent trajectories on the left side exit left, this
numerically verifies that the boundaries we have selected form an isolating
neighborhood.  Of course, this depends on the chosen grid, but by choosing a very fine
grid, we expect that this procedure provides a strong argument that we have found valid
boundaries. Some of the results from the integration of the tangent trajectories
are illustrated in 
Figure \ref{crtbpneighborhoodtest}.
Here, the trajectories are integrated from the boundary until they exit a
region a little larger than the boundary for easier visualization.  It can be
seen that the outer tangent trajectories exit the region quickly, while those
on the inner boundary wander into the isolating neighborhood before exiting
left.  In each case, it is verified that all tangent trajectories exit on the
appropriate side.

For the autonomous CRTBP, this procedure is sufficient to check the boundaries,
but for the non-autonomous ERTBP, additional computations are required.  In the
non-autonomous case, different tangent trajectories will be computed based on
which time, or epoch, is chosen.  The Jacobi constant is also no longer
available to simplify the computations.  Given these differences, the
boundaries must be checked using tangent trajectories at a range of time spans,
and a range of energies (or velocities) must also be checked.  Fortunately, for
the ERTBP, the system has a dimensionless period of $2\pi$, so it is only
necessary to check the boundaries on the interval from 0 to $2\pi$ before the
configuration repeats.  It is more difficult to determine the range of
velocities or energies to check, but an instantaneous Jacobi constant 
may be used as a guide.  This instantaneous Jacobi constant, and the equivalent
Hill's region, is computed from the instantaneous state in the pulsating ERTBP frame
when it is transferred to the CRTBP frame.  As the state varies in the ERTBP, the equivalent
Jacobi constant and Hill's region may be computed in the CRTBP, and these parameters
will vary as the infinitesimal mass moves along the trajectory.
This instantaneous Jacobi constant is only an
approximate guide, so during subsequent calculations, trajectories crossing the
isolating neighborhood boundaries are checked to verify that this assumed
envelope of velocities is not violated.

We can define the tangency computations more precisely using
the outer boundary as an example in the planar problem. 
The neighborhood boundary $\partial_R N$ can
be viewed in polar coordinates as $(r, \theta,s, \phi, \nu)$ where
$r$ is the radius of the cylinder, $\theta$ gives the angular location on
the cylinder, $s$ is the speed, $\phi$ gives the angular direction of the velocity,
and $\nu$ is the true anomaly.  Here, the right exit boundary is again
defined by the condition $r = 0.35$. The tangency condition is the condition
$\phi = \theta \pm \pi /2$.  A point $( .35, \ \theta_0,\  s, \phi_0,\  \nu)$
in the tangency set has two free variables; speed and true anomaly. For the
CRTBP the speed is determined by the Jacobi constant, and it may be calculated
as
\begin{equation}
s^2= 2 \Phi(r,\theta)-2C
.
\end{equation}

For the ERTBP, the Jacobi constant is also not constant but varies slowly, and a Jacobi function may be defined. The speed
of a solution at an exit point will be close to the speed of a solution
to the CRTBP at the same exit position.  The Jacobi function and the associated
Hill's zero velocity curve restricts the spacecraft positions in the rotating
coordinate system for a fixed Jacobi constant. If an interval $[C_0, C_1]$ of Jacobi
constants is used, then the smallest constant (or highest energy) determines a
barrier to the positions of a solution starting from an initial condition with $J = (C_0 + C_1)/2$
and the speed satisfies the constraint 
\begin{equation}
2 \Phi(r,\theta))-2C_0 \le s^2\le  2 \Phi(r,\theta)-2C_1
.
\end{equation}
The angle $\theta$ is also constrained by the inequality $0 \le \Phi (x,y) -
C_1$. Here, a Jacobi function in polar coordinates may be defined as
\begin{equation}
J(t, r, \theta, s, \phi) = \frac{1}{2} s^2  - \Phi (r, \theta)
.
\end{equation}
This function can provide some insight into how the Jacobi constant varies, but
we  have found that it is generally sufficient to use the instantaneous Jacobi
constant for our numerical computations.

An additional complication is the choice of boundaries to use relative to the
CRTBP boundaries in each ERTBP model.  In some of the models, the positions of
the primary and secondary can vary quite noticeably, making the choice of
boundary locations problematic.  From this perspective, the non-uniformly
rotating, pulsating ERTBP model provides the easiest gateway into studying this
problem.  In this model, the primary and secondary remain fixed, and we may
start our analysis by using fixed boundaries similar to the approach we used in
the CRTBP.

An initial test using the CRTBP boundaries 
in the non-uniformly rotating, pulsating ERTBP 
revealed that the inner boundary
failed the isolating neighborhood test when some of the tangent trajectories on
this boundary exited right.  By moving the inner boundary out to a radius of
1.02 centered on the barycenter, this problem was resolved for a range of
energies or velocities.  In this case, the boundaries for the ERTBP may now be
defined as follows.  The inner boundary $\partial_L N$ is defined by $r_L^2 =
x^2 + y^2$ where $r_L = 1.02$, and the outer boundary $\partial_R N$ is defined
as before by $r_R^2 = (x -1 + \mu)^2 + y^2$ where $r_R = 0.35$.
%
The range of required energies is defined iteratively using the instantaneous
CRTBP Jacobi constant as a convenient way to parameterize energy, keeping in
mind that the Jacobi constant does not exist in the ERTBP.  Subsequent
calculations of trajectories tracking orbits around L$_2$ are made, and the
energy, or Jacobi constant, of trajectories intersecting the isolating
neighborhood boundaries are computed.  It was determined from this process
that a range of Jacobi constants from 3.10 to 3.17 will be sufficient for
testing the isolating neighborhood boundaries.



\begin{figure}[ht!]
\centering
\subfigure[$\nu_i = 0$] {\label{exitnu0}
\includegraphics[width=0.48\textwidth]{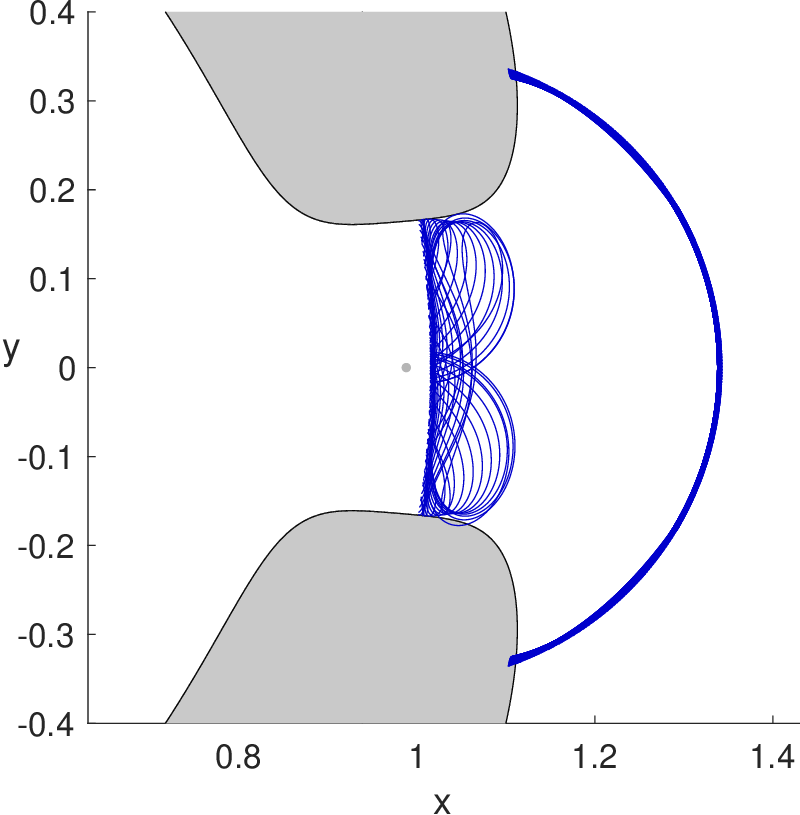}
}
\subfigure[$\nu_i = \frac{\pi}{2}$] {\label{exitnuhalfpi}
\includegraphics[width=0.48\textwidth]{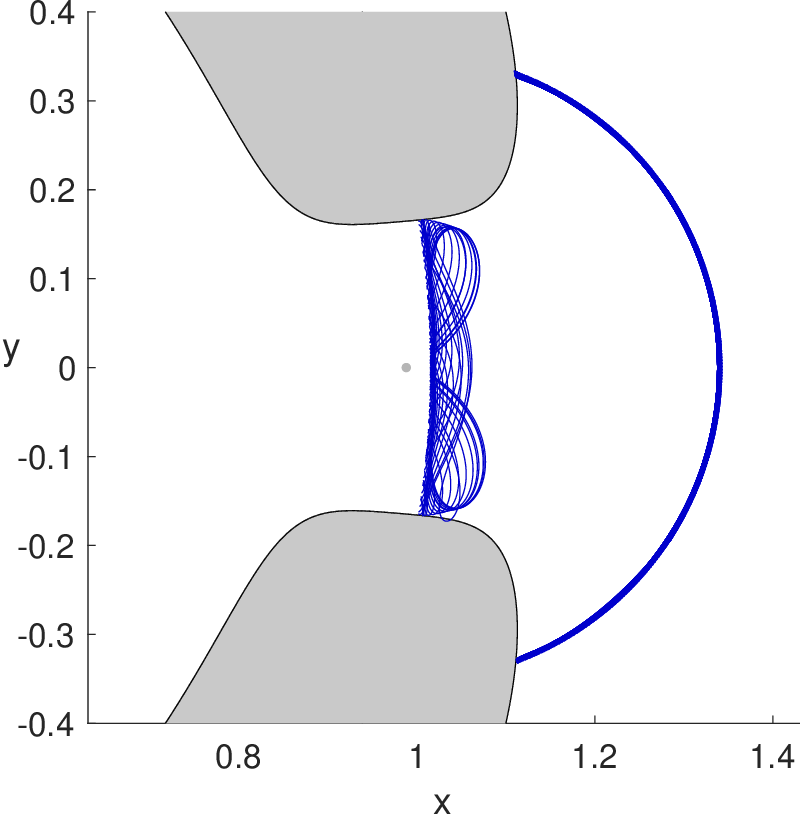}
}
\subfigure[$\nu_i = \pi$] {\label{exitnupi}
\includegraphics[width=0.48\textwidth]{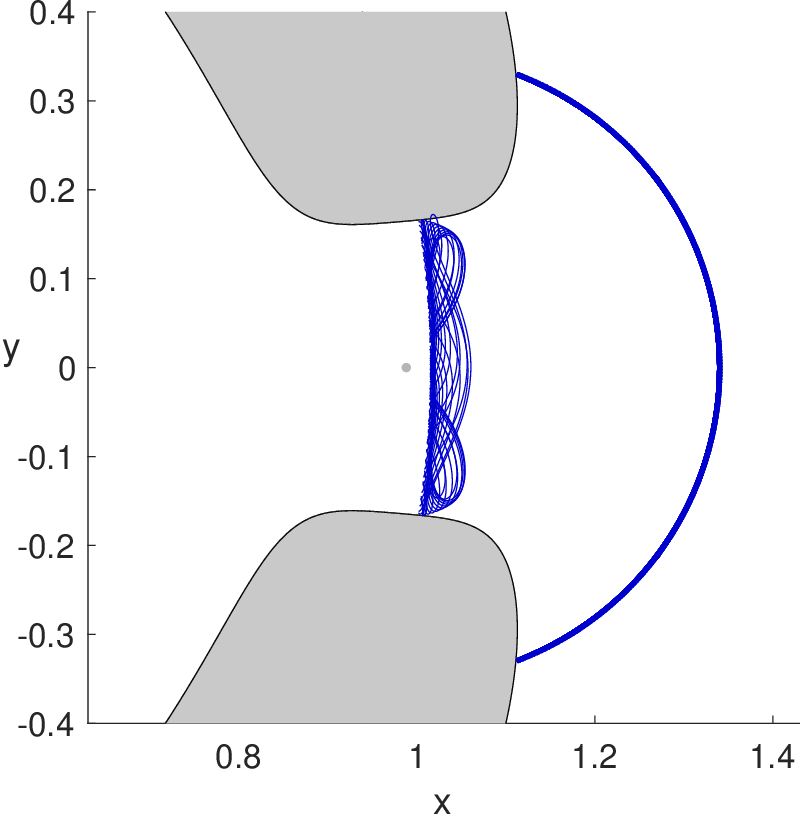}
}
\subfigure[$\nu_i = \frac{3\pi}{2}$] {\label{exitnu3halvespi}
\includegraphics[width=0.48\textwidth]{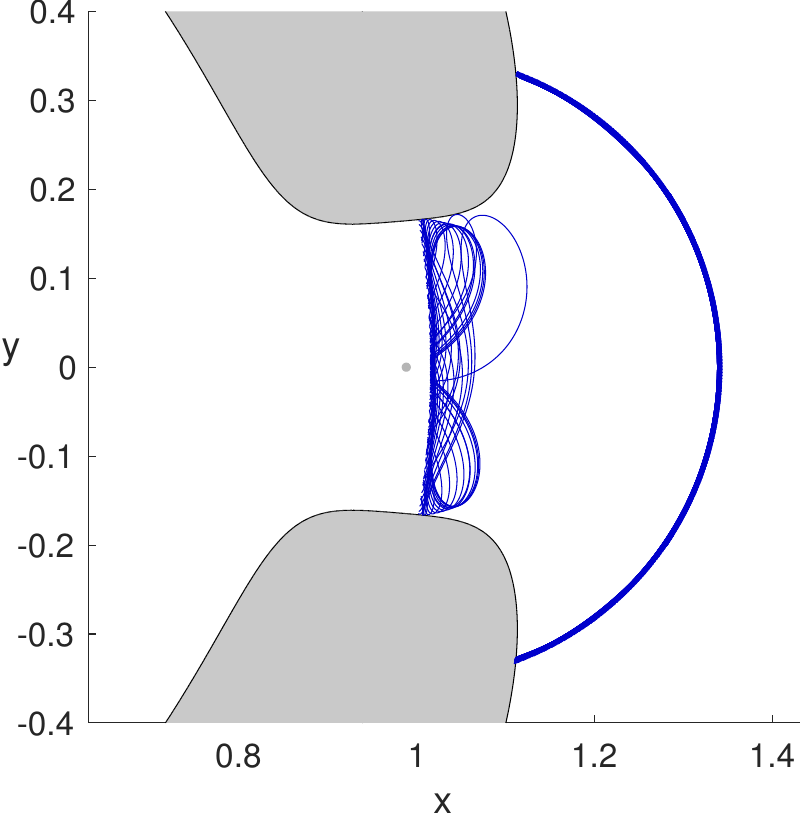}
}
\caption{
Test showing the tangent trajectories exiting the isolating neighborhood in the ERTBP.
(The Hill's regions for the CRTBP are plotted here for comparison only.)
}
\label{ertbpneighborhoodtest}
\end{figure}

\begin{figure}[ht!]
\centering
\subfigure[$C = 3.10$] {\label{exitnu0}
\includegraphics[width=0.48\textwidth]{exitnu0}
}
\subfigure[$C = 3.14$] {\label{c314boundarytestnu0}
\includegraphics[width=0.48\textwidth]{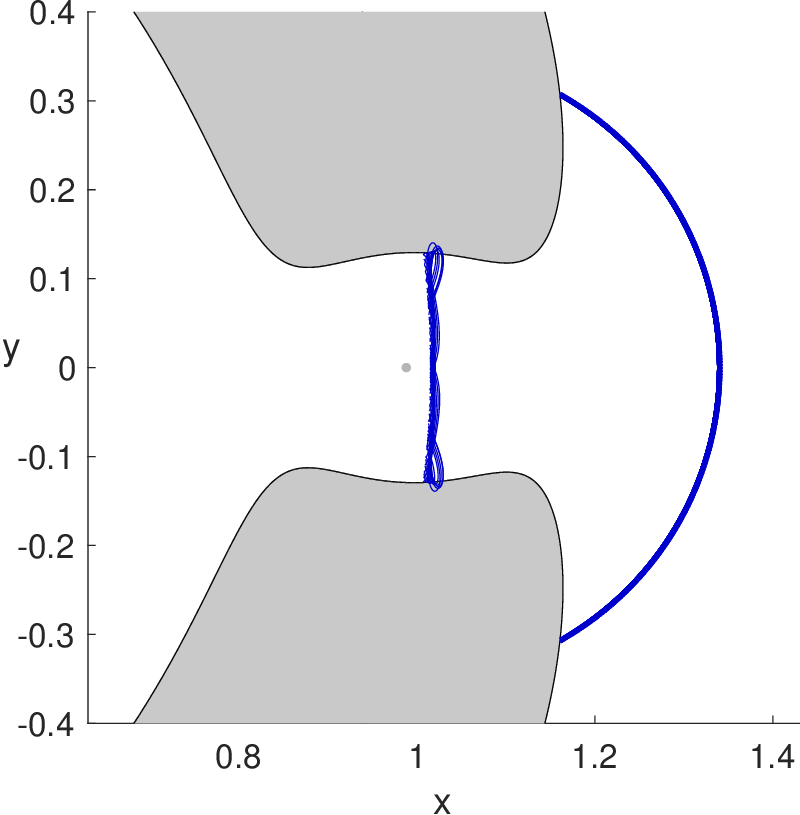}
}
\subfigure[$C = 3.17$] {\label{c314boundarytestnu0}
\includegraphics[width=0.48\textwidth]{c314boundarytestnu0}
}
\caption{
Test showing the tangent trajectories exiting the isolating neighborhood in the ERTBP for a range of energies with
$\nu_i = 0$.
(The Hill's regions for the CRTBP are plotted here for comparison only.)
}
\label{ertbpneighborhoodtestcrange}
\end{figure}

In order to verify that the boundaries we have selected act as isolating
neighborhood boundaries under all conditions required for the orbit tracking
procedure, the boundaries are tested for a range of initial dimensionless times
from 0 to $2\pi$ and a range of energies, parameterized using the instantaneous
CRTBP Jacobi constant, from C = 3.10 to 3.17.  In each case, it is checked that
the tangent trajectories on each boundary exit on that same boundary.  The most
stringent test is at the higher energy, corresponding to C = 3.10, and the
results for this case using four different initial true anomalies are shown in
Figure \ref{ertbpneighborhoodtest}.
Note that the instantaneous Hill's regions are provided for reference only since they do not
exist in the ERTBP, and as expected, some ERTBP trajectories cross over into
these regions.  It can be seen from this figure that some of the trajectories
on the left boundary occasionally travel into the interior before they exit
left, but all tested trajectories meet the criteria.  At lower energies, the
trajectories generally exit from the isolating neighborhood more quickly.  This
result can be seen in Figure \ref{ertbpneighborhoodtestcrange} where the tangent trajectories are shown
for several energies for the $\nu_i = 0$ case.

\clearpage

\section{Orbit Tracking Computations in the Planar ERTBP}

\begin{figure}[ht!]
\centering
\subfigure[Right-Left Exit Bisection] {\label{crtbprightleftexit}
\includegraphics[width=0.47\textwidth]{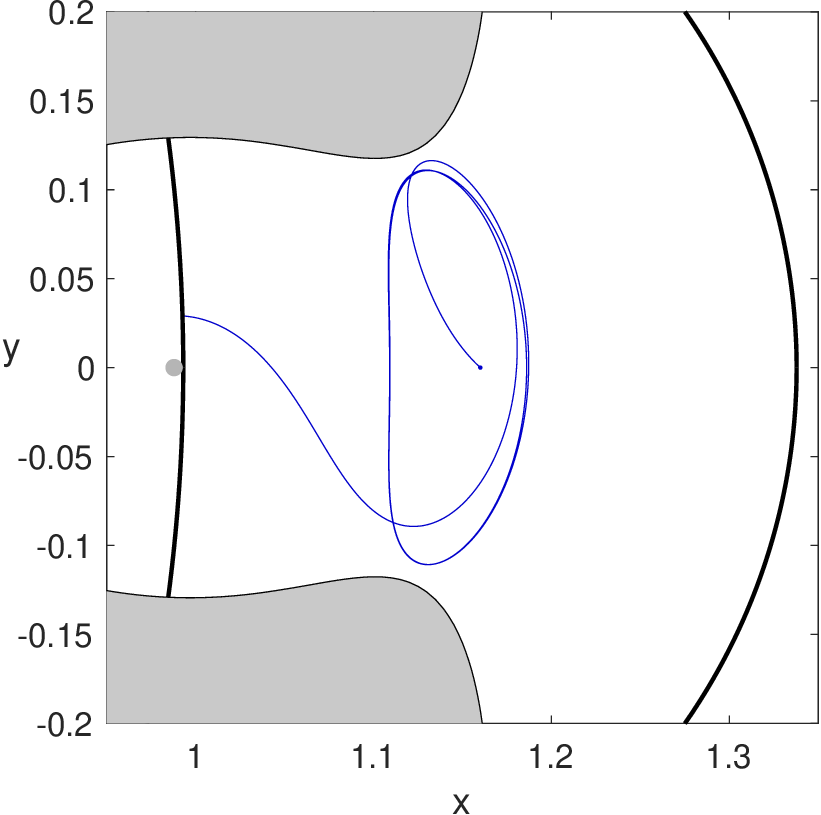}
}
\subfigure[Periodic Orbit From Orbit Tracking] {\label{crtbplyapunov}
\includegraphics[width=0.47\textwidth]{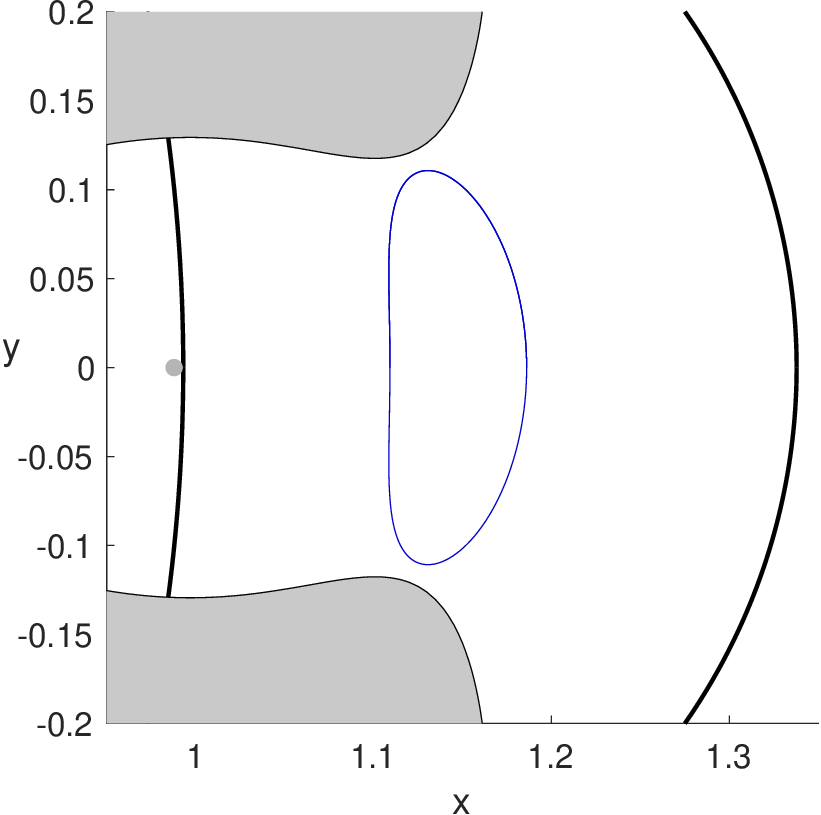}
}
\caption{
CRTBP forward asymptotic trajectory and quasiperiodic orbit computed for a base point
at (1.16, 0) with C = 3.14.
}
\label{crtbpplanarall}
\end{figure}

In our previous work, we used a bisection method to compute forward asymptotic
trajectories of the L$_2$ isolated invariant set in the CRTBP and then closely track
periodic and quasiperiodic orbits.\cite{Anderson:2019b}  In the planar CRTBP, this produces a
periodic orbit (the planar Lyapunov orbit), but in the ERTBP, we expect to
compute quasiperiodic orbits using this method.  As a first step, we compare
with the CRTBP by selecting a base point in the CRTBP at a chosen $(x,y)$ in the
computed isolating neighborhood at C = 3.14.  Various methods may be used as a
basis for the bisection algorithm, but for this case, we compute a circle of
velocities at the base point and find where they exit from the isolating
neighborhood.  The boundary between the sets of left and right exit
trajectories in Figure \ref{crtbpplanarall} may be computed using bisection for the base point at
(1.16, 0.0).
Two possible trajectories may be
obtained which, in this case, lie on the forward asymptotic set of the L$_2$
isolated invariant set.  Choosing one of these trajectories, and then computing
corrections using the orbit tracking algorithm at each $y = 0$ crossing\cite{Anderson:2017b}
produces the trajectory in Figure \ref{crtbplyapunov} which closely tracks the L$_2$
Lyapunov orbit.

Now a similar procedure may be performed in the non-uniformly rotating,
pulsating ERTBP model.  For this case, the same base point used in the CRTBP is
chosen, although it is now in a different frame.  The same velocities computed
in the CRTBP for C = 3.14 are also used at this point in the ERTBP, and a
similar procedure is performed to find both the left and right exit sets along
with the boundary between them.  Next, one of the forward asymptotic
trajectories on this boundary is selected, and an orbit tracking procedure is
performed at the same $y = 0$ crossings as in the CRTBP.  A slight modification
is required for the ERTBP however.  In the CRTBP, the Jacobi constant is used
to compute the velocity at each correction, insuring that the Jacobi constant
varies very little from the initial Jacobi constant.  Since there is no Jacobi
constant in the ERTBP, the speed that the trajectory intersects
$\Sigma_{y=0}$ with is used to compute the new corrected velocity using bisection.
%
%
%
%
%
\begin{figure}[ht!]
\centering
\subfigure[Right-Left Exit Bisection] {\label{ertbpc314bisect}
\includegraphics[width=0.47\textwidth]{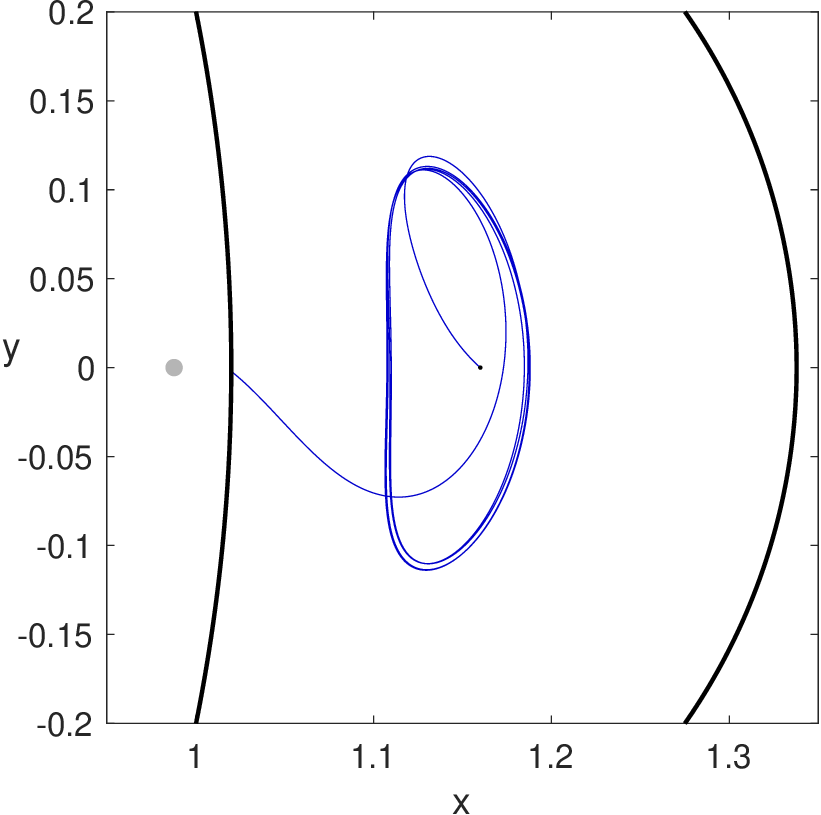}
}
\subfigure[Quasiperiodic Orbit From Orbit Tracking] {\label{ertbpc314orbit}
\includegraphics[width=0.47\textwidth]{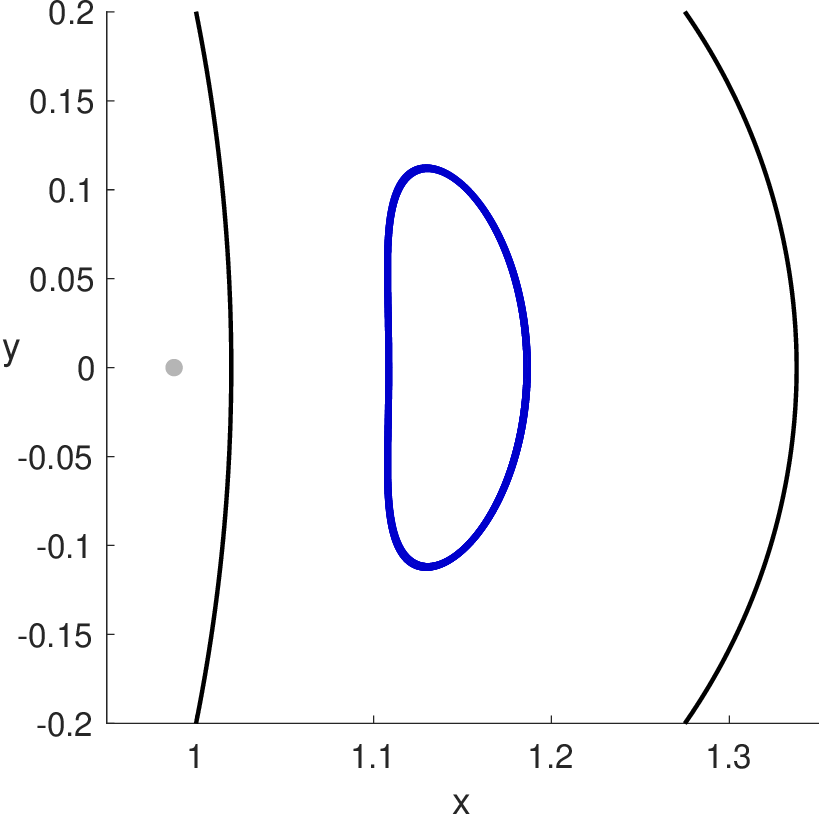}
}
\caption{
ERTBP forward asymptotic trajectory and quasiperiodic orbit computed for a base point
at (1.16, 0) with C = 3.14.
}
\label{ertbpplanarall}
\end{figure}
The results from this procedure are shown in Figure \ref{ertbpplanarall}.  
Note that the initial approach from
the base point is shown in Figure \ref{ertbpc314bisect}, and only the trajectories closely tracking the quasiperiodic orbit
are plotted in Figure \ref{ertbpc314orbit}.
In
this case, the computed trajectory in Figure \ref{ertbpc314orbit} approaches a
quasiperiodic orbit, or a 2-torus, rather than the periodic Lyapunov orbit
obtained in the CRTBP.  These results are in line with what would be expected
in the planar ERTBP from Capinski and Zgliczy\'nski,\cite{Capinski:2018} and it
is similar to results found in the bicircular problem.\cite{Jorba:1997, Jorba:2020, Rosales:2021}
It is worth mentioning that some of these results require that the perturbation be sufficiently small, but our
method does not specifically require this.
The characteristics of the quasiperiodic orbit may be seen more clearly by plotting
a Poincar\'e section using $\Sigma_{y=0}$ as shown in Figure \ref{ertbpc314xxdot}.
\begin{figure}[ht!]
\centering
\includegraphics[width=0.47\textwidth]{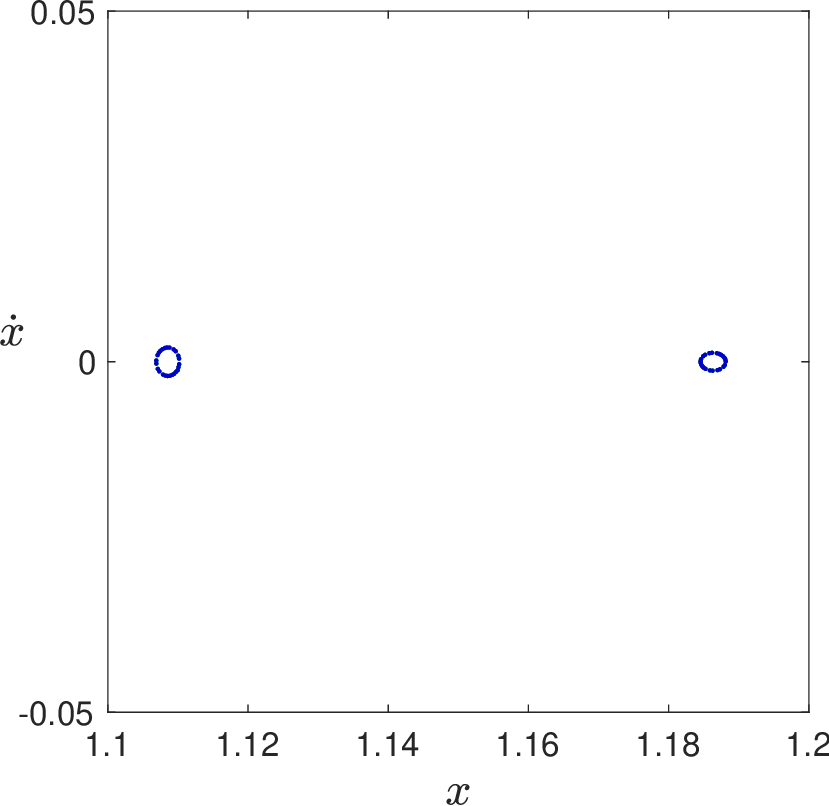}
\caption{
Intersections of the quasiperiodic orbit in Figure \ref{ertbpc314orbit} with the surface of section $\Sigma_{y=0}$.
}
\label{ertbpc314xxdot}
\end{figure}
Now that we have implemented the algorithm to track the quasiperiodic orbit,
we next check the trajectories used in the computation to ensure that they do not
go outside the range of validity of the isolating neighborhood boundaries that we
have already tested.  To test this, we take all of the trajectories used in the bisection
portion of the orbit tracking algorithm and compute their energies, parameterized
by the Jacobi constant, as they intersect the boundaries.  Performing this computation
gives $C_{min} \approx 3.1077$ and $C_{max} \approx 3.1675$ which is inside
the range of energies when the boundaries were checked earlier.

\begin{figure}[ht!]
\centering
\subfigure[$\nu_i = \pi/2$] {\label{orbithalfpi}
\includegraphics[width=0.47\textwidth]{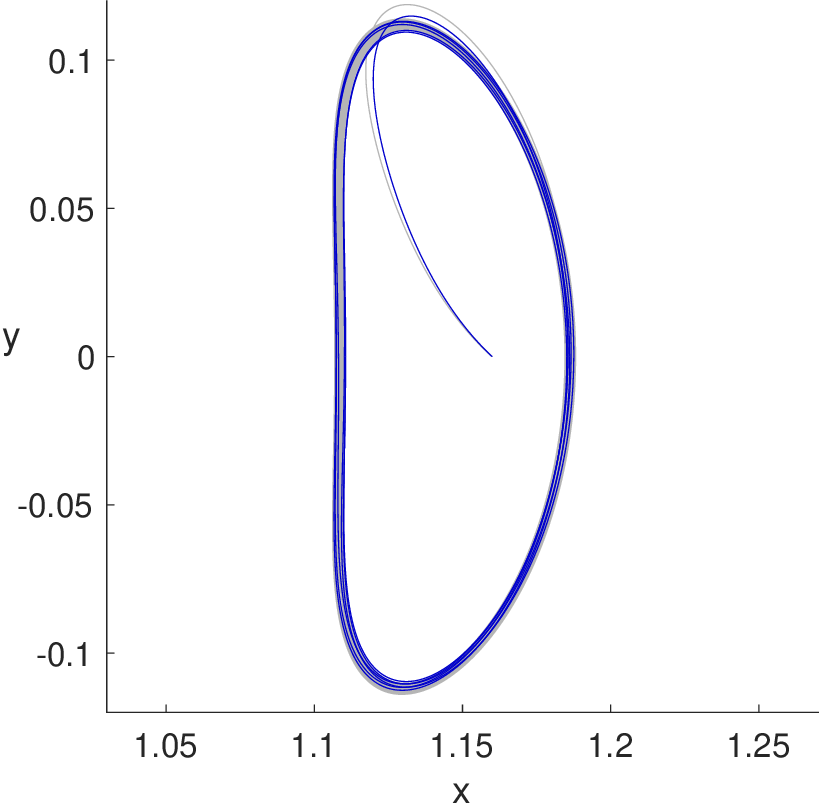}
}
\subfigure[$\nu_i = \pi$] {\label{orbitpi}
\includegraphics[width=0.47\textwidth]{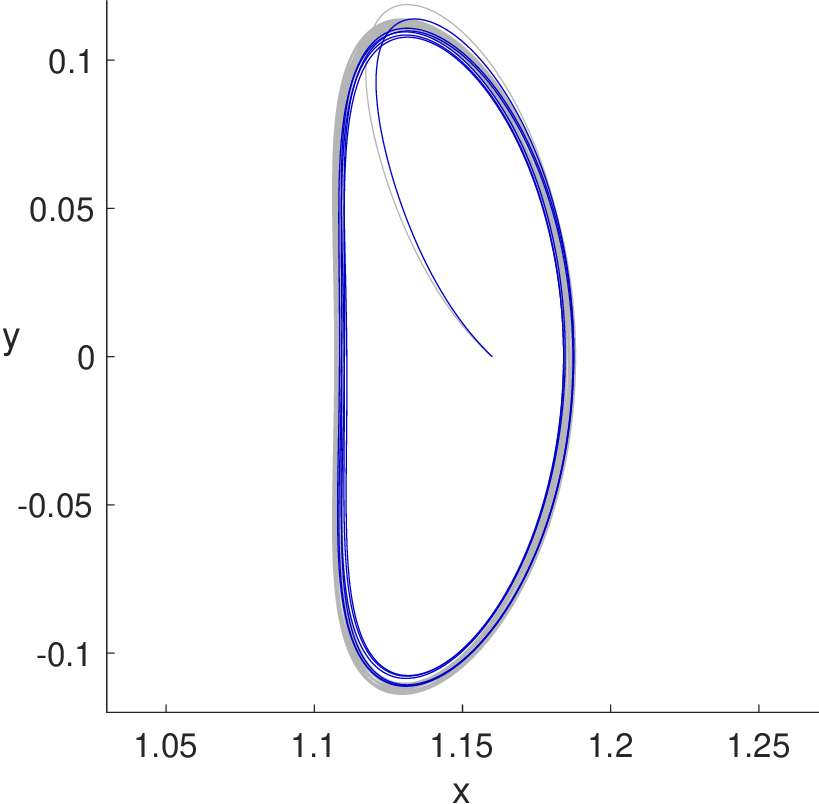}
}
\caption{
Quasiperiodic orbit for base point at (1.16, 0) and C = 3.14 compute for different initial
$\nu_i$ values.  The gray orbit in the background was computed at $\nu_i = 0$, and it is
provided for reference.
}
\label{configothernus}
\end{figure}

One additional factor to include in the computation of orbits in the ERTBP is
the initial epoch of the trajectory.  The orbit in Figure \ref{ertbpplanarall}
was computed for an initial epoch with $\nu_i = 0$.  If a different initial
epoch is chosen, a different quasiperiodic orbit will be computed.  Orbits for
$\nu_i = \pi/2, \pi$ are shown in Figure \ref{configothernus}.
In each case, the orbit computed at $\nu_i = 0$ is plotted in gray in the background
for reference.  The largest apparent difference is found for $\nu_i = \pi$.




\clearpage
\section{Orbit Tracking Computations in the Spatial ERTBP}

\begin{figure}[ht!]
\centering
\subfigure[Inner Cylinder] {\label{innercylinder}
\includegraphics[width=0.47\textwidth]{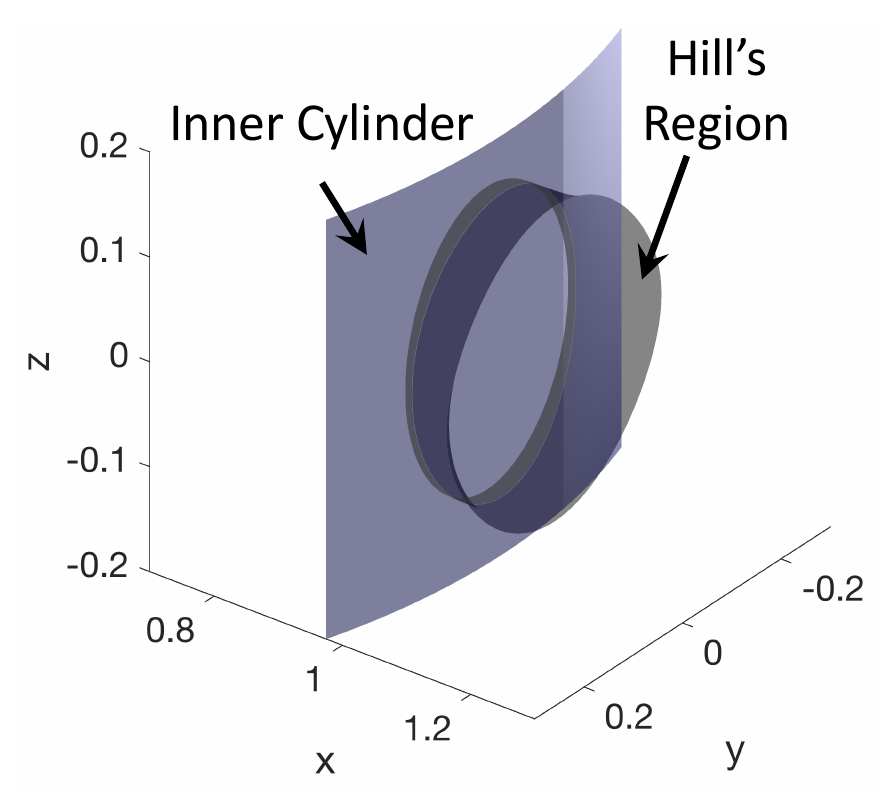}
}
\subfigure[Outer Cylinder] {\label{zerocylc310}
\includegraphics[width=0.37\textwidth]{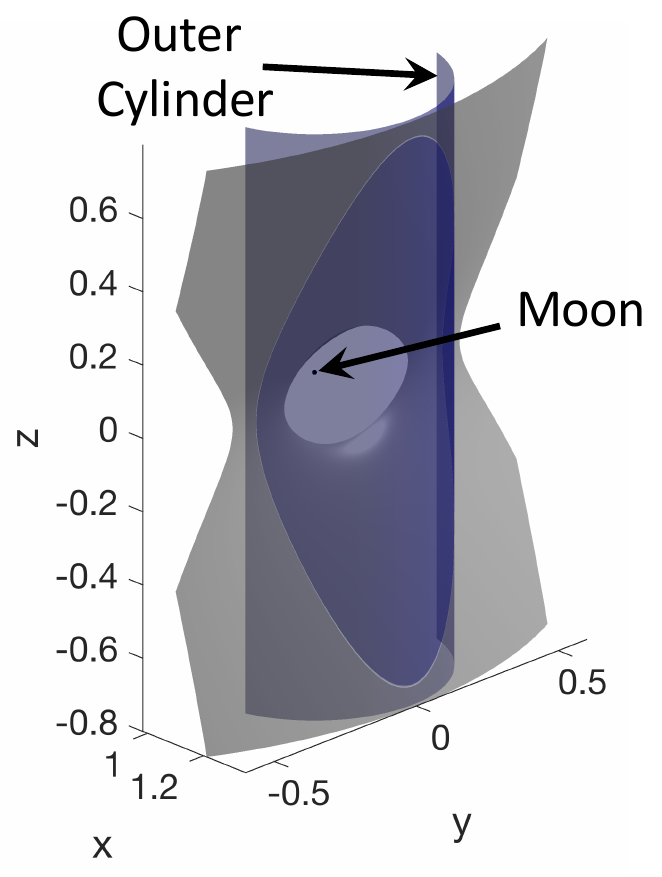}
}
\caption{
Schematics showing the inner and outer isolating neighborhood cylinder boundariers relative
to the Hill's region for $C = 3.10$.  Note that the Hill's region is computed 
in the CRTBP and is shown for convenience to visualize the location of 
the isolating neighborhood boundaries.
}
\label{allzerocylc310}
\end{figure}

The methods implemented for the planar ERTBP may be extended to the spatial
ERTBP by first using cylinders rather than circles for the isolating
neighborhood boundaries.  The boundaries may be checked using methods very
similar to those implemented in Anderson, Easton, and Lo.\cite{Anderson:2019b}
In brief, the checks first require the computation of a grid of points on the
outer and inner cylinders.  Tangent trajectories are integrated forward and
backward at each of the grid points to ensure that the tangent trajectories
exit on the side they started on.  For the ERTBP, these conditions were checked
for a range of initial times between 0 and 2$\pi$ and for $3.10 \le C \le 3.17$
where $C$ is once again used here as a convenient way to parameterize the
velocities.  The outer and inner cylinders used for $C = 3.10$ are shown relative to the Hill's
region 
computed from the instantaneous Jacobi constant in the CRTBP
in Figure \ref{allzerocylc310}.
In this case, it is important to verify that the cylinder intersects the Hill's region
at the top and the bottom.  


For the spatial ERTBP  the seven-dimensional phase space with coordinates 
$(\nu, x,y,z,x^\prime, y^\prime, z^\prime )$ and
the differential Equations \ref{eomx} through \ref{eomz} are used for the numerical computations.
A Jacobi manifold $M(C)$ is the set of all points in phase space satisfying Equation
\ref{jacobiconstant}. The Jacobi manifolds form a layer in phase space: 
\begin{equation}
\Lambda = \bigcup \{ M(C) : 3.10 \le C \le 3.17 \}
.
\end{equation}
The neighborhood $N$ in phase space in which we search for quasiperiodic
solutions is defined by three conditions:
\begin{equation}
r_L^2 = x^2 +y^2 \ge 1.02^2
\end{equation}
\begin{equation}
r_R^2 = (x+1 -\mu)^2 +y^2 \le .35^2
,
\end{equation}
and 
\begin{equation}
(\nu, x,y,z,x^\prime, y^\prime, z^\prime ) \in \Lambda
.
\end{equation}
As will be described next, trial solutions are generated by choosing initial
conditions inside $N$  which belong to the Jacobi manifold $M(3.14)$. In our
numerical experiments we observe that these solutions do not exit the layer
$\Lambda$ before they exit the neighborhood $N$.

Once the isolating neighborhood boundaries have been checked, specific base points 
for computing trial solutions within
the isolating neighborhood may be examined.  We explore several base points here
in the $z = 0$ surface of section $(\Sigma_{z=0})$, focusing on the cases where $\nu_i = 0$.
We know from previous analyses in the CRTBP, that if the base point is in a region
where periodic or quasiperiodic orbits exist, we can find points 
that stay in the isolating neighborhood forward and backward in time.  The procedure
for finding these points is detailed in Anderson, Easton, and Lo,\cite{Anderson:2020b}
but a brief summary is given here.  As a first step, spheres of velocities are computed
for a selected basepoint, and the velocities on these spheres that exit left and right both forward and backward in time are first computed.  The boundaries of the left and right exit trajectories are computed, and their intersections contain the states that remain in the isolating neighborhood
both forward and backward in time.
In the ERTBP, the velocity magnitude used for the velocity spheres is chosen based
upon a selected instantaneous Jacobi constant.
These points are then states for trajectories contained
within the isolated invariant set, and we may use them as initial conditions
for tracking periodic or quasiperiodic orbits within the isolated invariant set.

\begin{figure}[ht!]
\centering
\includegraphics[width=0.47\textwidth]{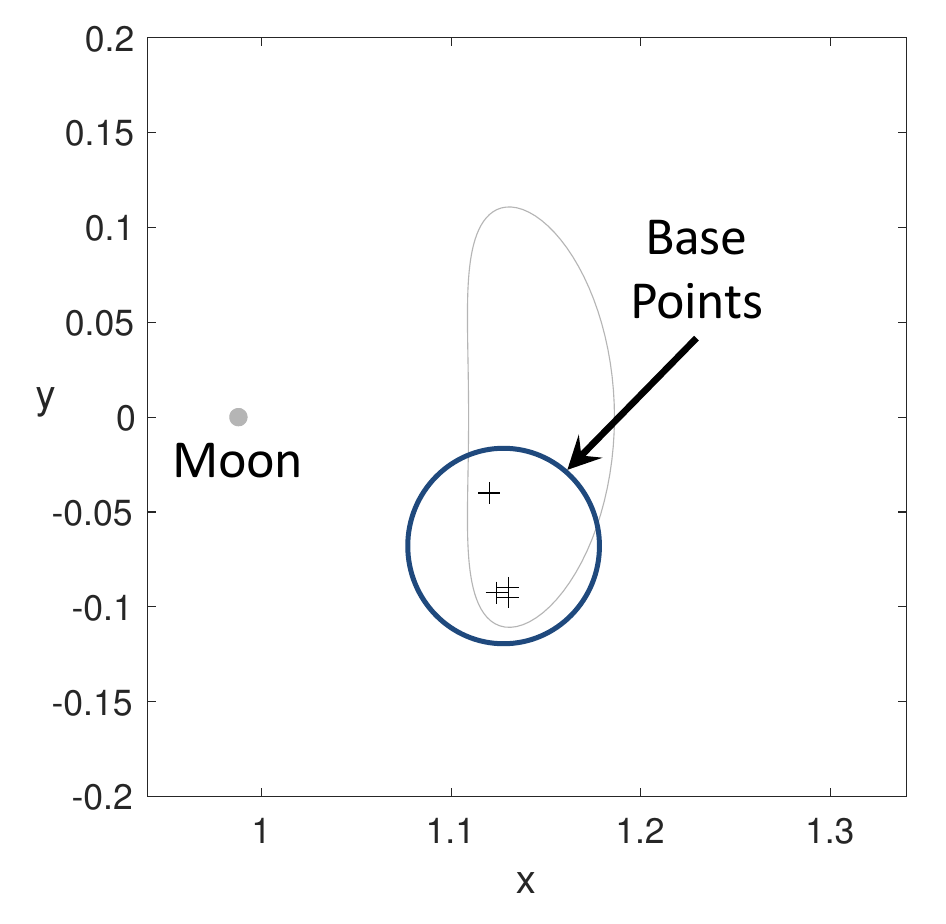}
\caption{
Overview of the base points used to track orbits in the ERTBP with C = 3.14.  (The gray
Lyapunov orbit is computed using the CRTBP and is provided for reference only.)
}
\label{spatialoverview}
\end{figure}

Using the initial conditions computed using the intersections on the velocity sphere
at a particular base point, several trajectories are integrated forward and corrected
to track orbits that stay in the region.  In each case, the correction is applied to the spatial
orbit as it intersects $\Sigma_{z=0}$, and it compensates for the limited precision of 
the initial conditions and numerical integration.  The correction in each case is applied
to the velocity, and it is the smallest  velocity correction necessary to move back
to the asymptotic set.
These corrections are generally very small, and specific correction magnitudes for
much longer orbits will be discussed later.
The selected base points are plotted in Figure \ref{spatialoverview} relative to the Lyapunov
orbit computed at C = 3.14 in the CRTBP (plotted only for reference since our
computations are in the ERTBP).  
Points across the region were selected so as 
to obtain different types of orbits corresponding to the Lissajous and quasi-halo
orbits in the CRTBP.

\begin{figure}[ht!]
\centering
\subfigure[Orbit in Configuration Space] {\label{c314casec}
\includegraphics[width=0.41\textwidth]{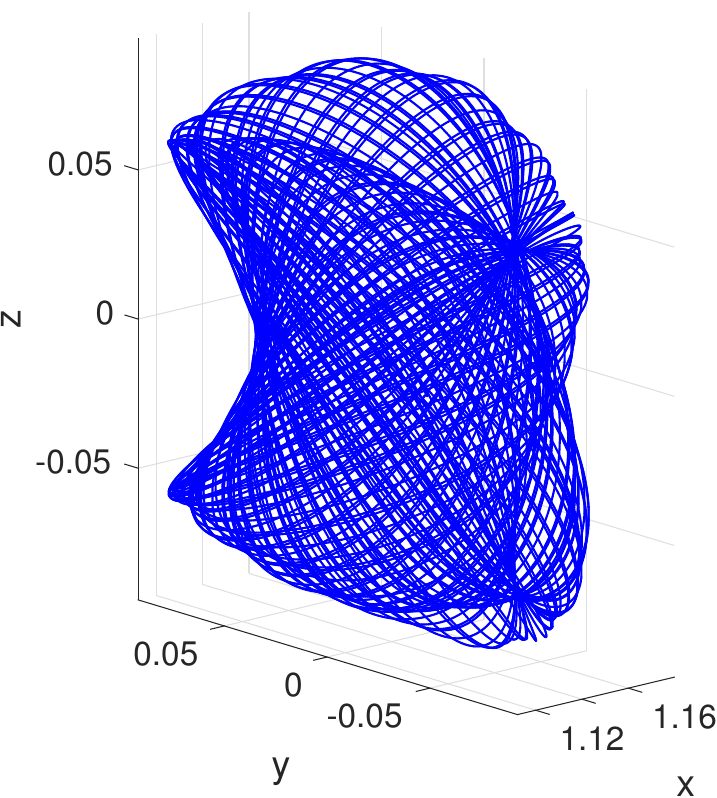}
}
\subfigure[$\Sigma_{z=0}$ Intersections] {\label{c314casec}
\includegraphics[width=0.50\textwidth]{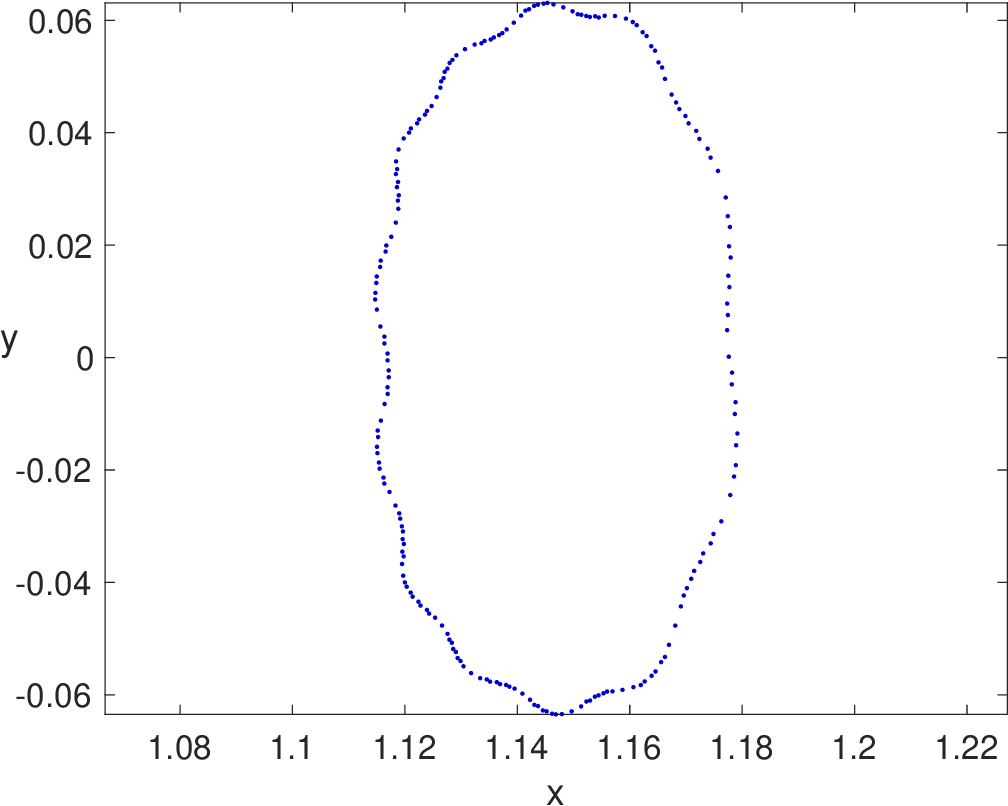}
}
\caption{
Quasiperiodic orbit tracked starting from base point $(1.12, -0.04)$ at C = 3.14.
}
\label{allc314casec}
\end{figure}

The case with a base point at $(1.12, -0.04)$ is shown in Figure
\ref{allc314casec} with $\nu_i = 0.0$.  (Unless otherwise specified all of the
orbit tracking results will start with $\nu_i = 0.0$.)
Here, the orbit is similar to the Lissajous orbits found in the CRTBP.
Examining the intersections with $\Sigma_{z=0}$ reveals a similar topology to
that found in the CRTBP, but now there are oscillations around a curve that
would normally be smooth.  
%
In this case, the tori generated in the ERTBP are three-dimensional tori rather
than the two-dimensional tori generated in the CRTBP.

\begin{figure}[ht!]
\centering
\subfigure[Orbit in Configuration Space] {\label{c314casea}
\includegraphics[width=0.47\textwidth]{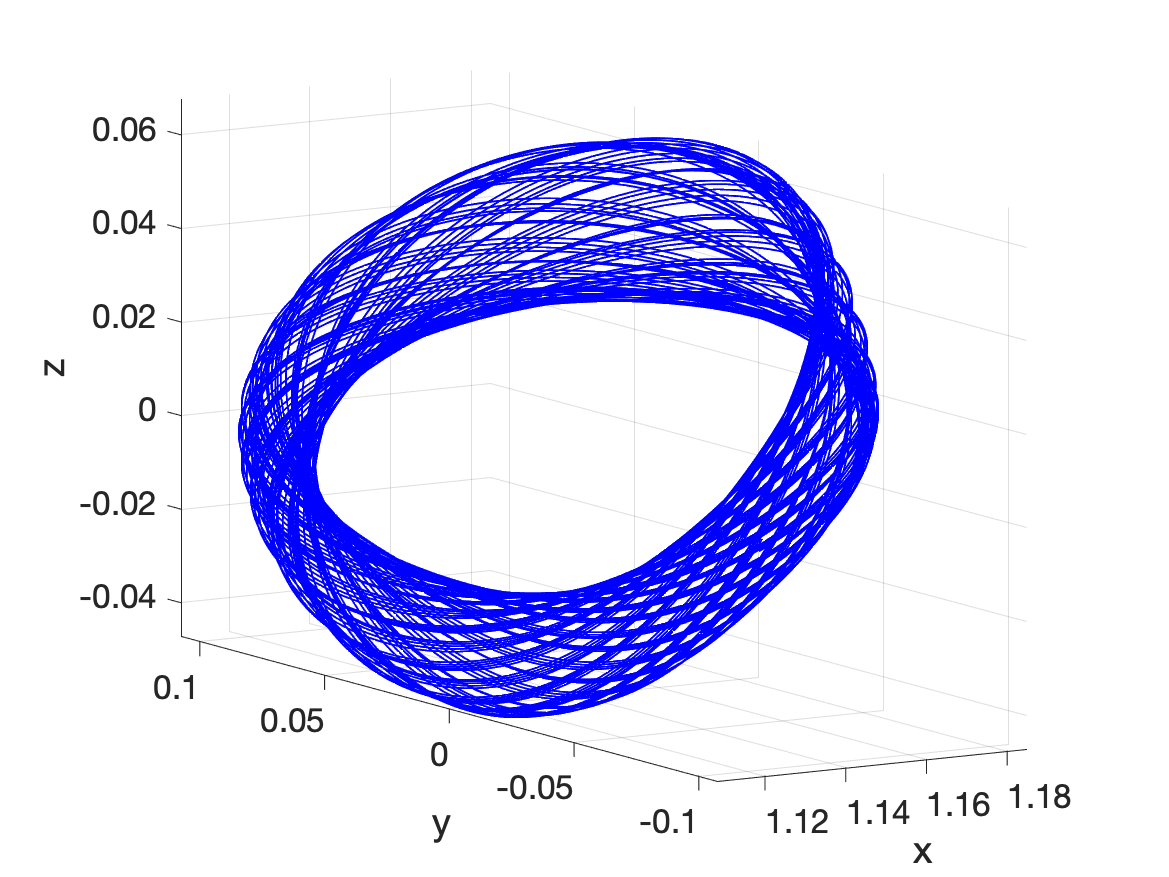}
}
\subfigure[$\Sigma_{z=0}$ Intersections] {\label{caseaxxdotc314}
\includegraphics[width=0.47\textwidth]{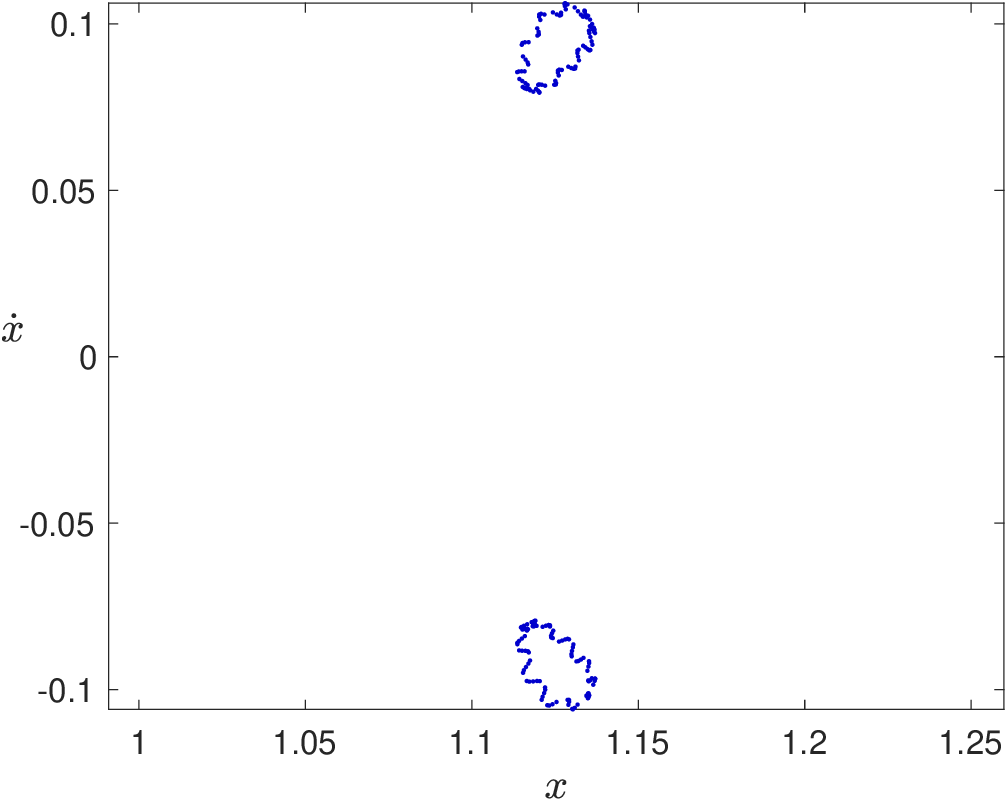}
}
\caption{
Quasiperiodic orbit tracked starting from base point $(1.13, -0.09)$ at C = 3.14.
}
\label{allc314casea}
\end{figure}

\begin{figure}[ht!]
\centering
\subfigure[Orbit in Configuration Space] {\label{c314caseb}
\includegraphics[width=0.47\textwidth]{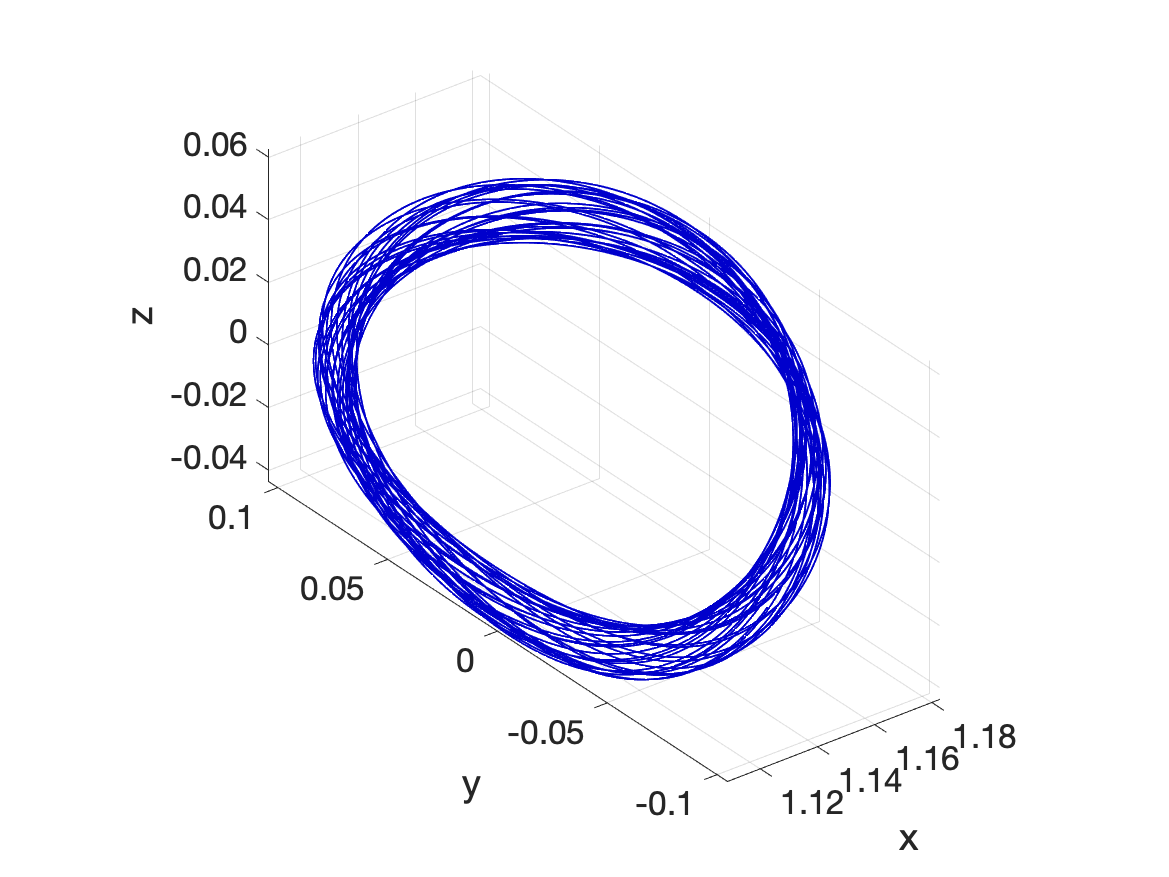}
}
\subfigure[$\Sigma_{z=0}$ Intersections] {\label{c314caseb}
\includegraphics[width=0.47\textwidth]{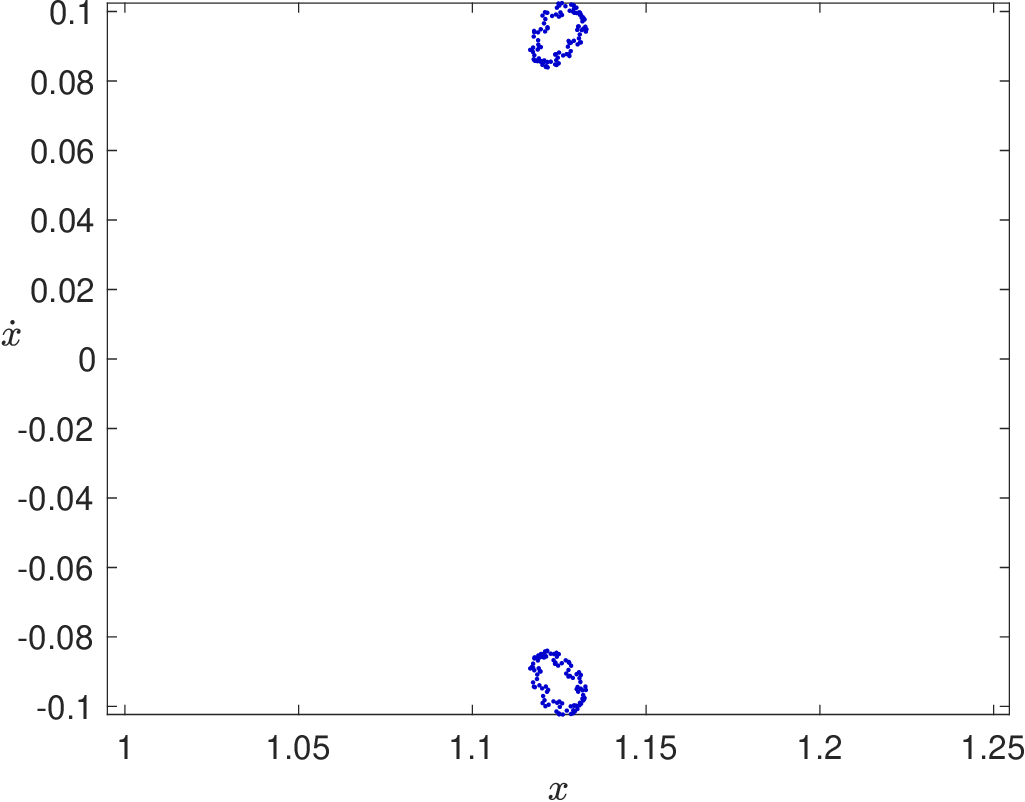}
}
\caption{
Quasiperiodic orbit tracked starting from base point $(1.13, -0.095)$ at C = 3.14.
}
\label{allc314caseb}
\end{figure}

\begin{figure}[ht!]
\centering
\subfigure[Orbit in Configuration Space] {\label{c314cased}
\includegraphics[width=0.47\textwidth]{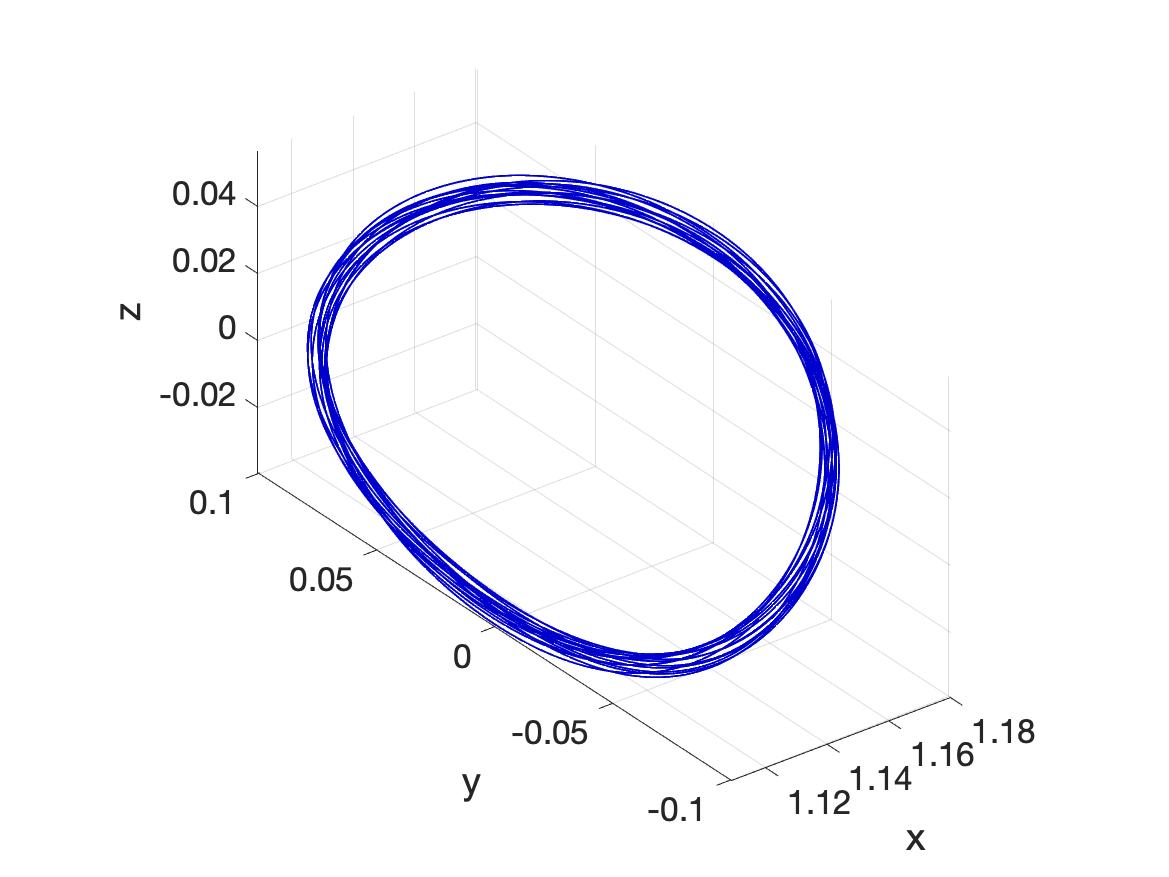}
}
\subfigure[$\Sigma_{z=0}$ Intersections] {\label{c314cased}
\includegraphics[width=0.47\textwidth]{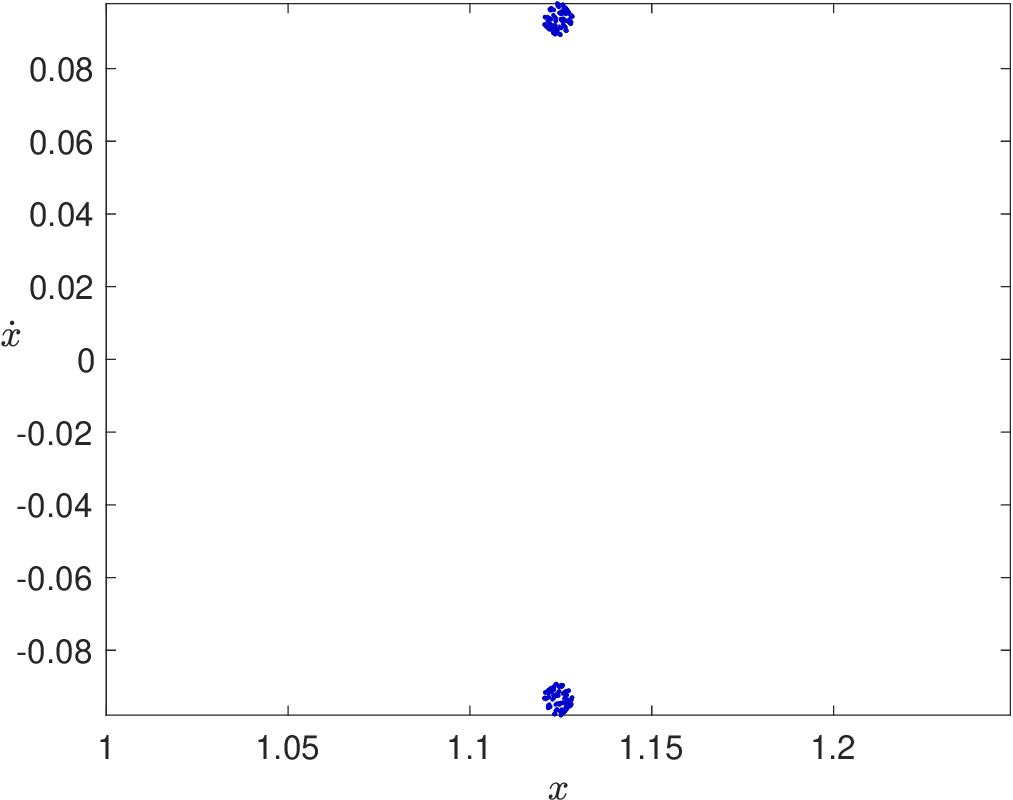}
}
\caption{
Quasiperiodic orbit tracked starting from base point $(1.124, -0.0925)$ at C = 3.14.
}
\label{allc314cased}
\end{figure}

Next, three base points were tracked in the quasi-halo region, and
the results are given in Figures \ref{allc314casea} through \ref{allc314cased}.
%
%
In each case, the first portions of the quasi-halo orbits are shown so as to be
able to still observe the structure in the orbits.  An initial set of
intersections with $\Sigma_{z=0}$ are also shown in the Poincar\'e sections.
These base points produce three quasi-halo orbits of different sizes, and in
each case the Poincar\'e sections show a similar overall structure to the ones
observed in the CRTBP.
In each case, there are also still oscillations about the expected curve.

\begin{figure}[ht!]
\centering
\subfigure[Orbit in Configuration Space] {\label{orbitx1.1300y-0.0900}
\includegraphics[width=0.47\textwidth]{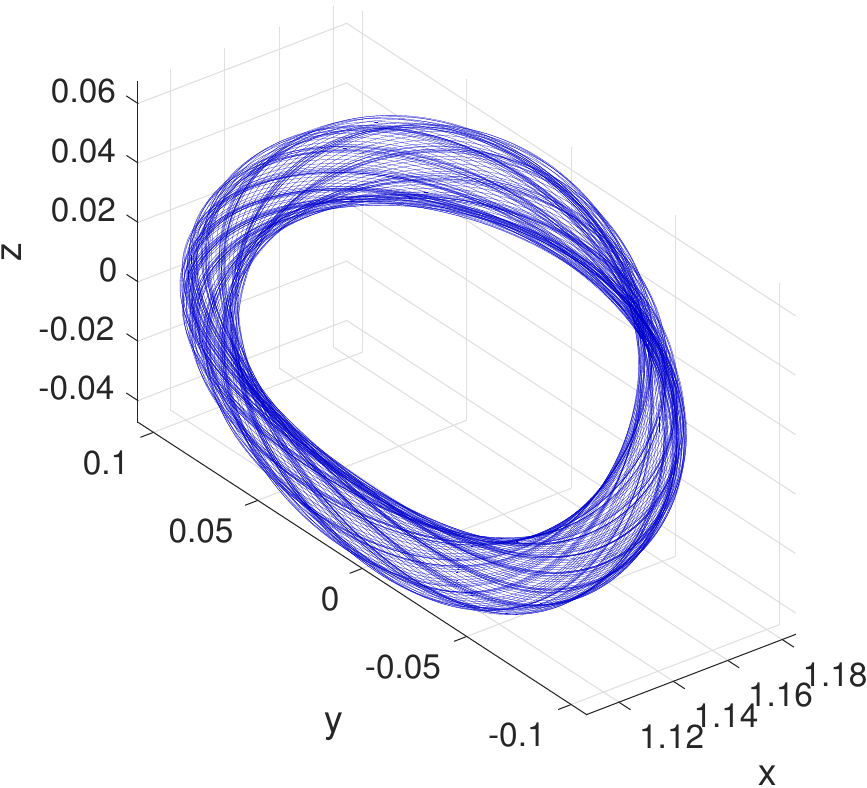}
}
\subfigure[$\Sigma_{z=0}$ Intersections] {\label{poinx1.1300y-0.0900}
\includegraphics[width=0.45\textwidth]{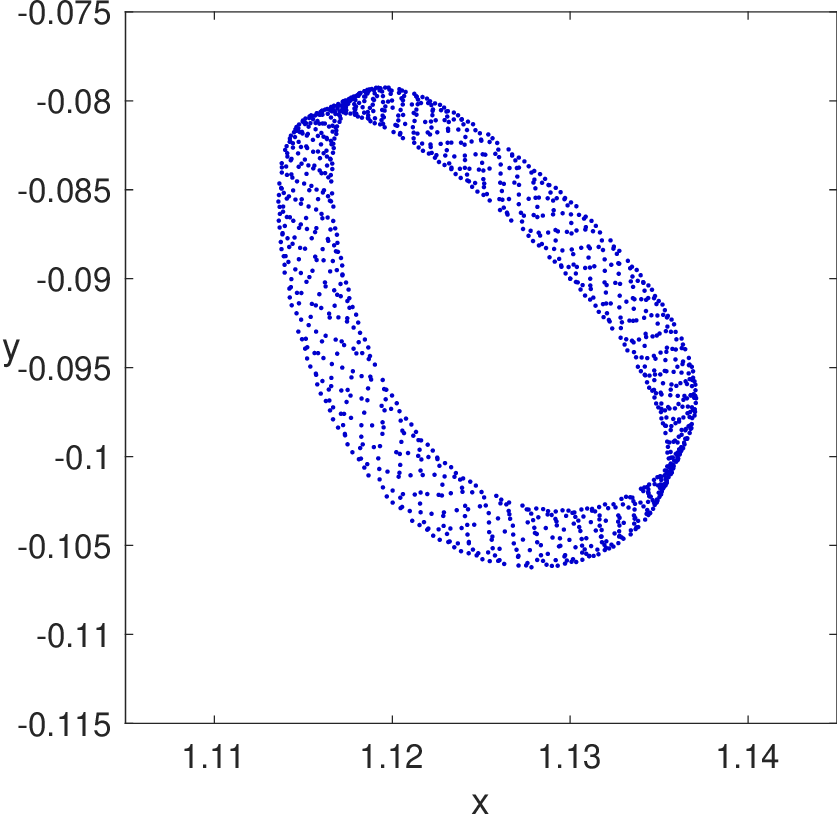}
}
\caption{
Quasiperiodic orbit tracked starting from base point $(1.130, -0.0900)$ with $\nu_i=0$ at C = 3.14
with additional points.
}
\label{longorbitx1.1300y-0.0900}
\end{figure}

\begin{figure}[ht!]
\centering
\subfigure[Orbit in Configuration Space] {\label{orbitx1.1300y-0.0902}
\includegraphics[width=0.47\textwidth]{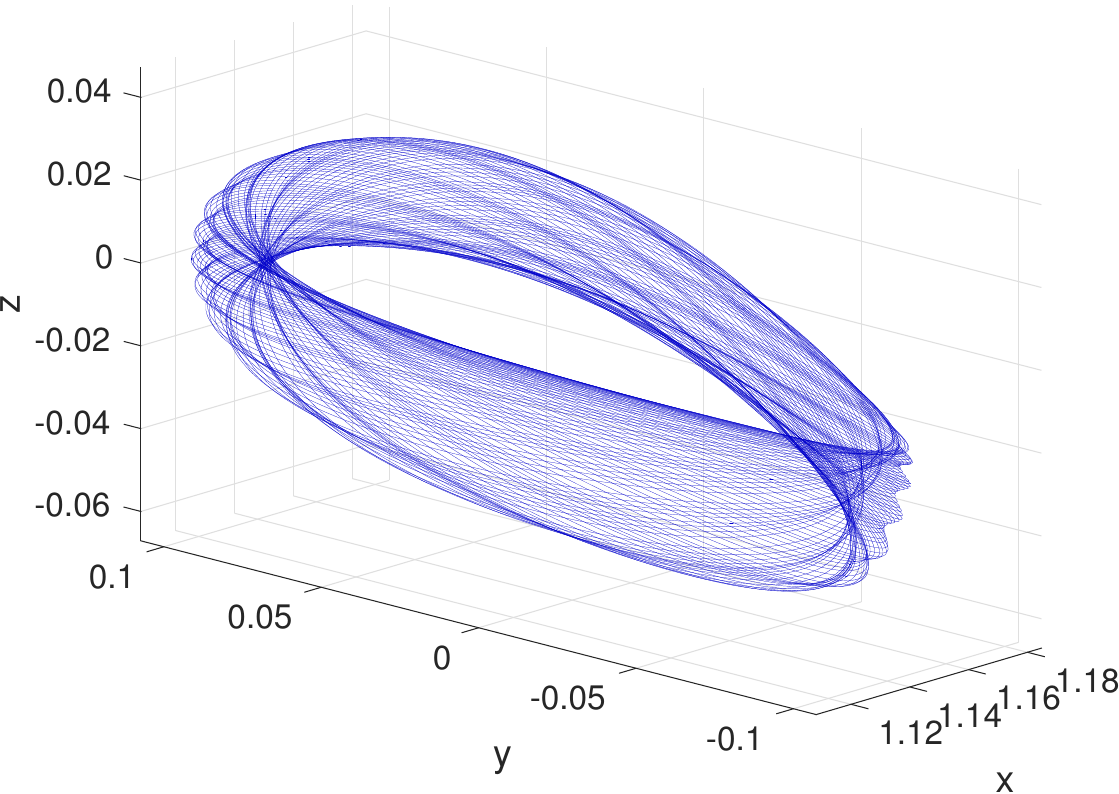}
}
\subfigure[$\Sigma_{z=0}$ Intersections] {\label{poinx1.1300y-0.0902}
\includegraphics[width=0.40\textwidth]{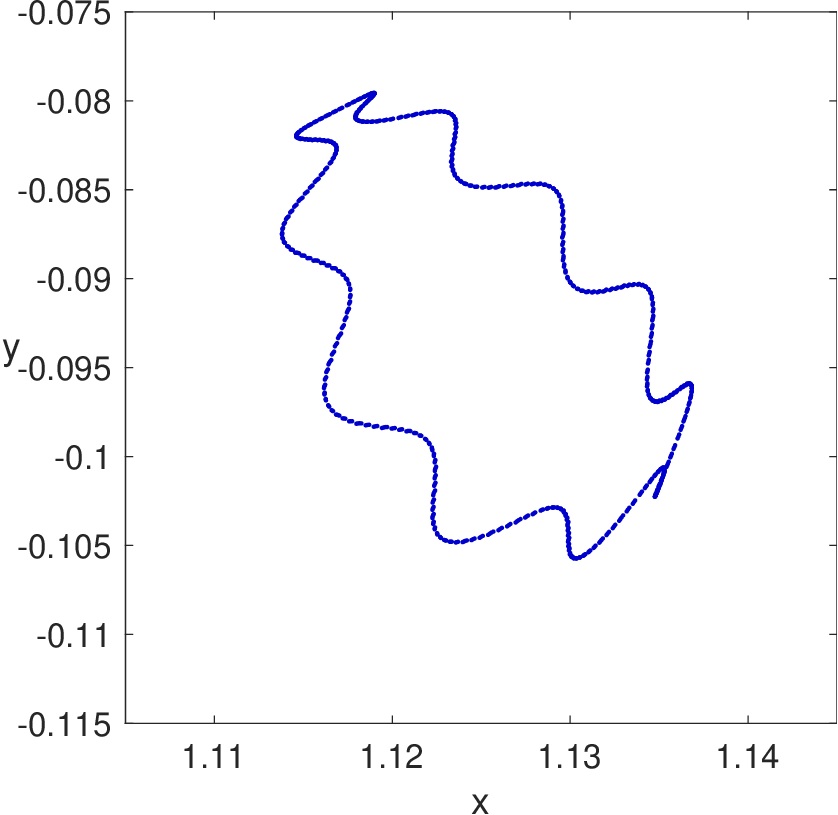}
}
\caption{
Quasiperiodic orbit tracked starting from base point $(1.130, -0.0902)$ at C = 3.14 with $\nu_i=0$.
}
\label{longorbitx1.1300y-0.0902}
\end{figure}

\begin{figure}[ht!]
\centering
\subfigure[Orbit in Configuration Space] {\label{orbitx1.1300y-0.0903}
\includegraphics[width=0.47\textwidth]{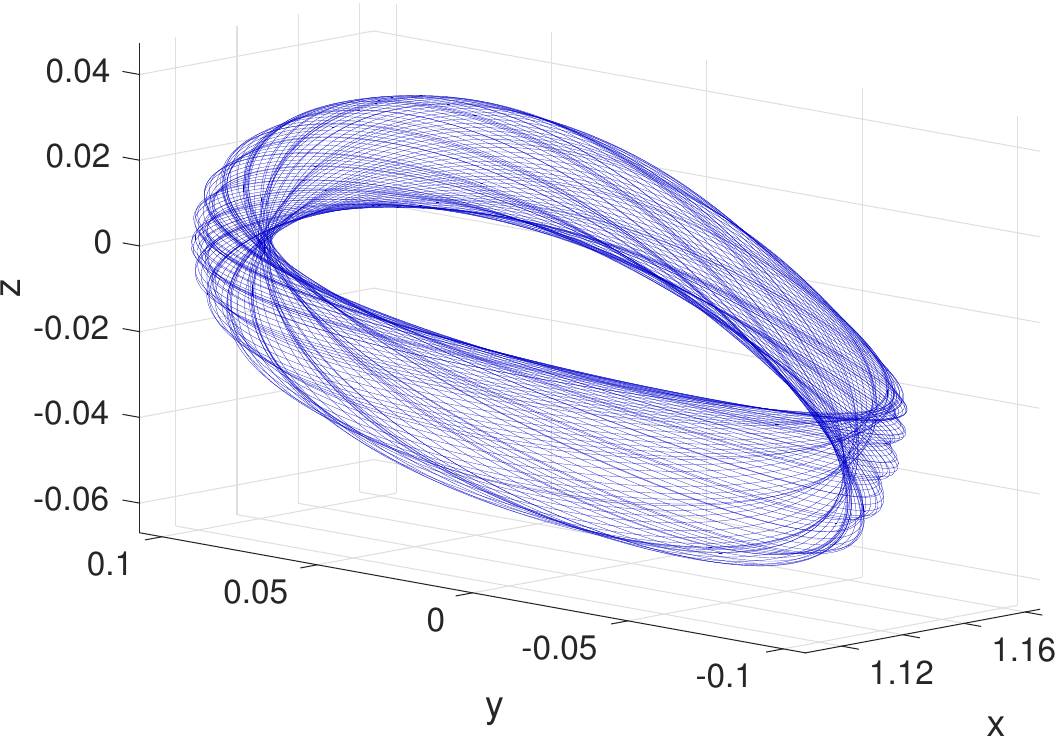}
}
\subfigure[$\Sigma_{z=0}$ Intersections] {\label{poinx1.1300y-0.0903}
\includegraphics[width=0.40\textwidth]{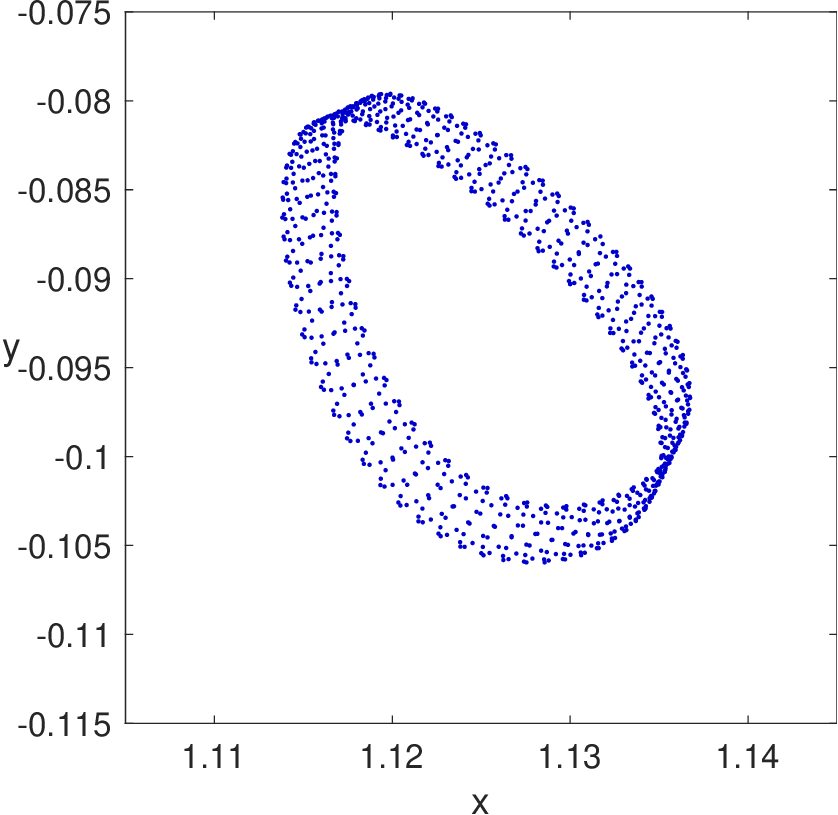}
}
\caption{
Quasiperiodic orbit tracked starting from base point $(1.130, -0.0903)$ at C = 3.14 with $\nu_i =0$.
}
\label{longorbitx1.1300y-0.0903}
\end{figure}

If the orbit tracking algorithm is continued further for additional
intersections with $\Sigma_{z=0}$, some additional interesting behavior may be
observed for different base points.  The difference in behavior between
different base points was explored, and three different results are shown for
base points that are relatively close in Figures \ref{longorbitx1.1300y-0.0900}
through \ref{longorbitx1.1300y-0.0903}.  Some characteristics of each
of these cases are compiled in Table \ref{orbitchar}.  
\begin{table}
\caption{
Orbit tracking characteristics for large revolution orbit calculations
}
\label{orbitchar}
\begin{center}
\begin{tabular}{c c c c} 
\hline
Base Point & Revolutions & Dimensionless Total $\Delta$V & Dimensional Total $\Delta$V (m/s)\\
\hline
(1.130, -0.0900)& 650 & $6.15\times10^{-5}$ & $\approx 0.063$\\
(1.130, -0.0902)& 500 & $4.68\times10^{-5}$ & $\approx 0.048$\\
(1.130, -0.0903)& 500 & $5.43\times10^{-5}$ & $\approx 0.056$\\
\hline
\end{tabular}
\end{center}
\end{table}
Even with over 500 revolutions, the positions are continuous,
and the sum of all velocity discontinuities from the orbit tracking algorithm was always less then $6.2 \times 10^{-5}$ or
appoximately 0.063 m/s.
In each case, the initial base point
was varied only slightly, but significantly different behavior was found.  Note that only
the portion of the Poincar\'e section with $x < 0$ is shown to more easily see the structure
in the intersecting points.
Figure \ref{longorbitx1.1300y-0.0900} uses the same initial base point as
Figure \ref{allc314casea}, but in this case significantly more revolutions were added.
There is initial oscillatory behavior that fills in as a band as the number of
revolutions is increased.  If the base point location is changed slightly to
(1.130, -0.0902) as shown in Figure \ref{longorbitx1.1300y-0.0902} then the points lie
on a curve, and more structure is seen in the orbit in configuration space.
Finally, if another nearby base point is examined at (1.130, -0.0902), as shown
in Figure \ref{longorbitx1.1300y-0.0903}, another type of oscillatory behavior
is found.
%
%

\begin{figure}[ht!]
\centering
\subfigure[Orbit in Configuration Space] {\label{diftimeorbitx1.1300y-0.0902}
\includegraphics[width=0.47\textwidth]{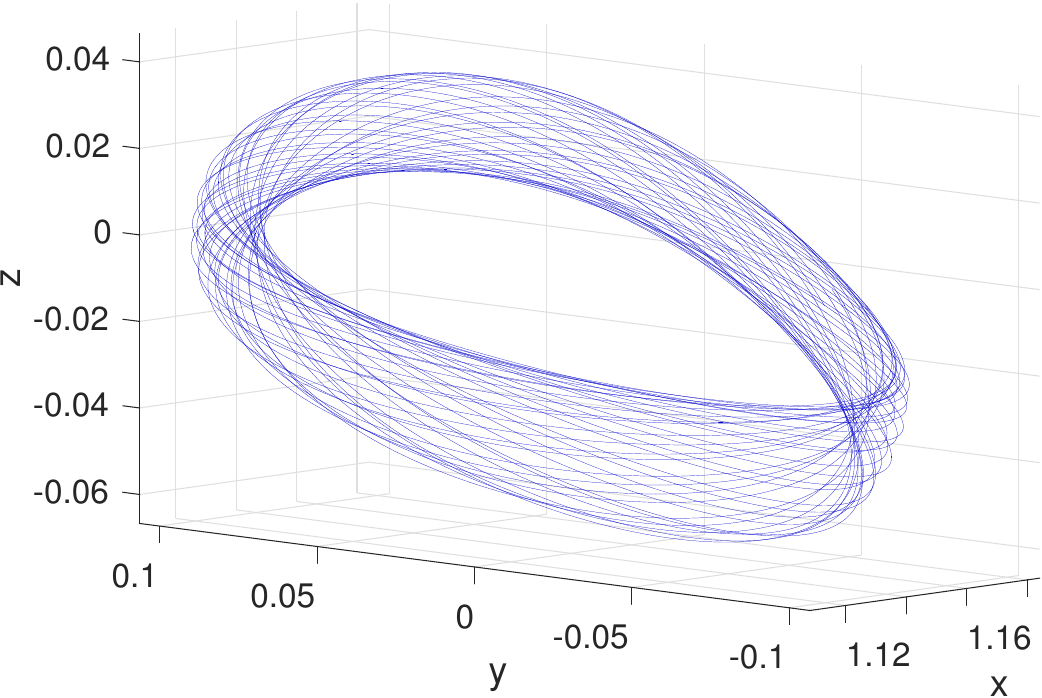}
}
\subfigure[$\Sigma_{z=0}$ Intersections] {\label{diftimepoinx1.1300y-0.0902}
\includegraphics[width=0.40\textwidth]{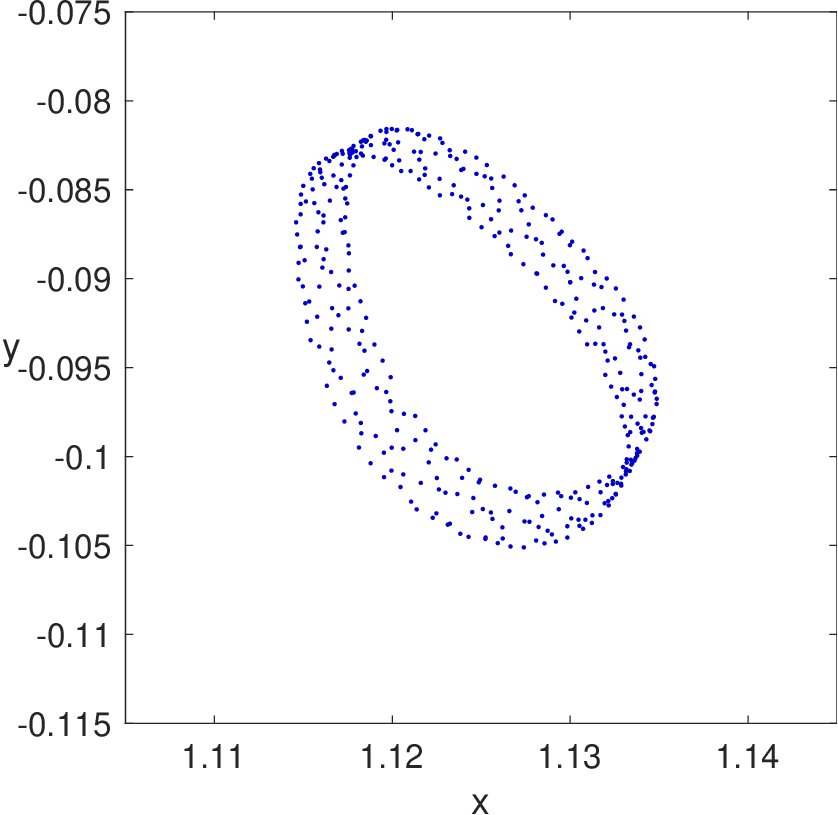}
}
\caption{
Quasiperiodic orbit tracked starting from base point $(1.130, -0.0902)$ at C = 3.14 with $\nu_i = \pi/2$.
}
\label{difftimelongorbitx1.1300y-0.0902}
\end{figure}

These results depend of course on the initial time that the orbit tracking algorithm
is initiated, and if a different initial time is selected, then different behavior
is observed.  The results for a base point at $(1.130, -0.0902)$ computed with an initial
time $\nu = \pi/2$ are shown in Figure \ref{difftimelongorbitx1.1300y-0.0902}.
While the case with $\nu_i = 0$ formed a line in $\Sigma_{z=0}$, with the new
initial time the intersections form a band.  This result is more similar to the
results using base point (1.130, -0.0900) with $\nu_i = 0$, and reinforces the
importance of the initial epoch in the computation of the resulting
quasiperiodic orbits.


\clearpage

\subsection{View of Spatial Orbits in Different Coordinate Frames}

\begin{figure}[ht!]
\centering
\subfigure[Orbit in Configuration Space] {\label{nonpulseconfigx1.1300y-0.0902}
\includegraphics[width=0.49\textwidth]{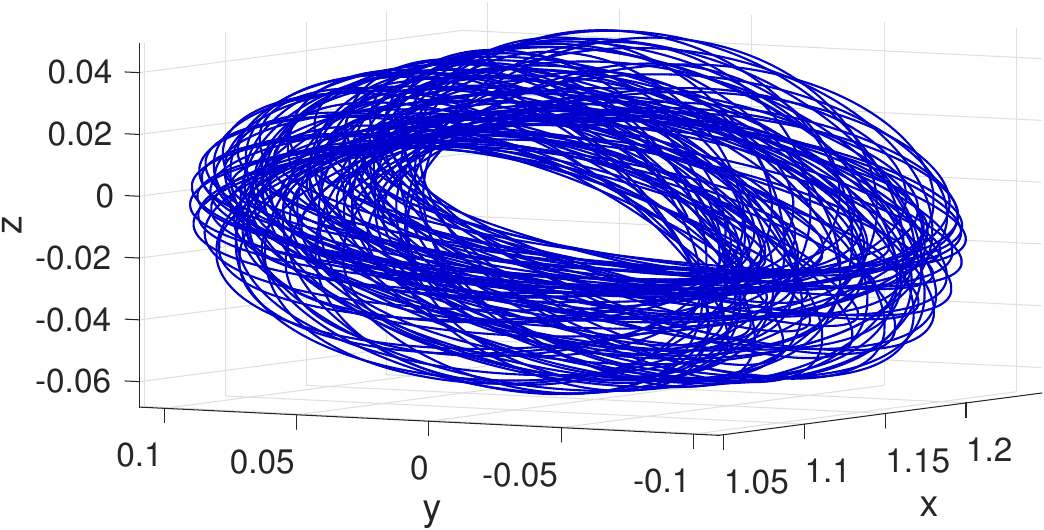}
}
\subfigure[$\Sigma_{z=0}$ Intersections] {\label{nonpulsexxdotx1.1300y-0.0902}
\includegraphics[width=0.40\textwidth]{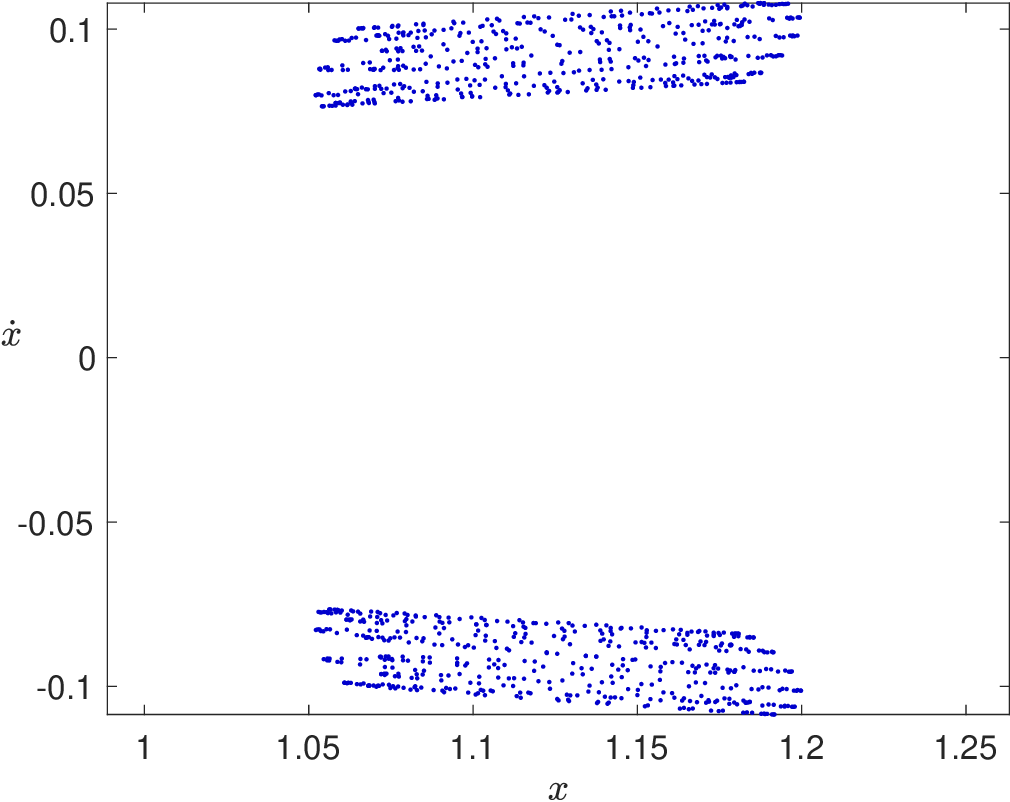}
}
\caption{
Quasiperiodic orbit tracked starting from base point $(1.13, -0.0902)$ at 
C = 3.14 in the variable rotating, non-pulsating frame.
}
\label{nonpulseall}
\end{figure}

In the pulsating frame, the orbits do not appear too dissimilar from the orbits
computed in the CRTBP.  In real-world computations, the frame is not pulsating,
and it is interesting to compute these orbits in different frames.  An orbit
computed in the non-pulsating, variable rotating frame corresponding to the
orbit in Figure \ref{allc314casea} is shown in Figure \ref{nonpulseall}.
Here, the orbit still has the general characteristics of a quasi-halo orbit,
but as might be expected, the orbit is elongated in the direction of the
$x$-axis.  This is easily seen in the Poincar\'e section shown in Figure
\ref{nonpulsexxdotx1.1300y-0.0902} where the curve in Figure \ref{poinx1.1300y-0.0902} is elongated
as well, and structure is less easily discerned.

\begin{figure}[ht!]
\centering
\subfigure[Orbit in Configuration Space] {\label{constantrotconfig}
\includegraphics[width=0.49\textwidth]{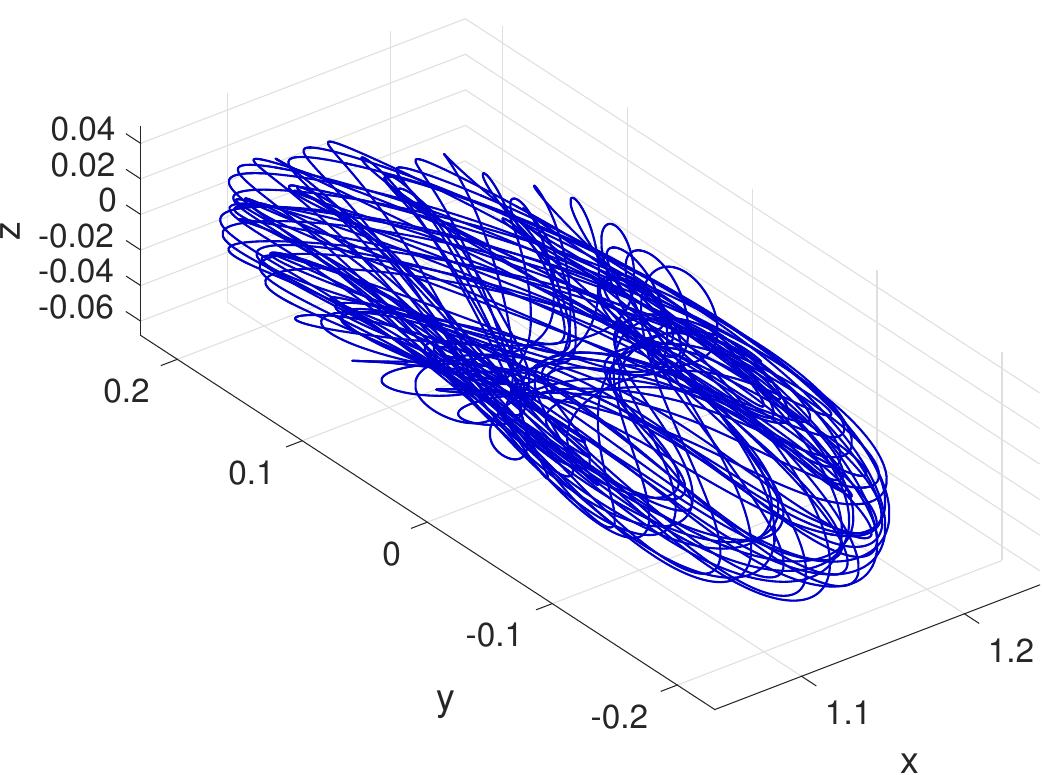}
}
\subfigure[$\Sigma_{z=0}$ Intersections] {\label{constantrot}
\includegraphics[width=0.40\textwidth]{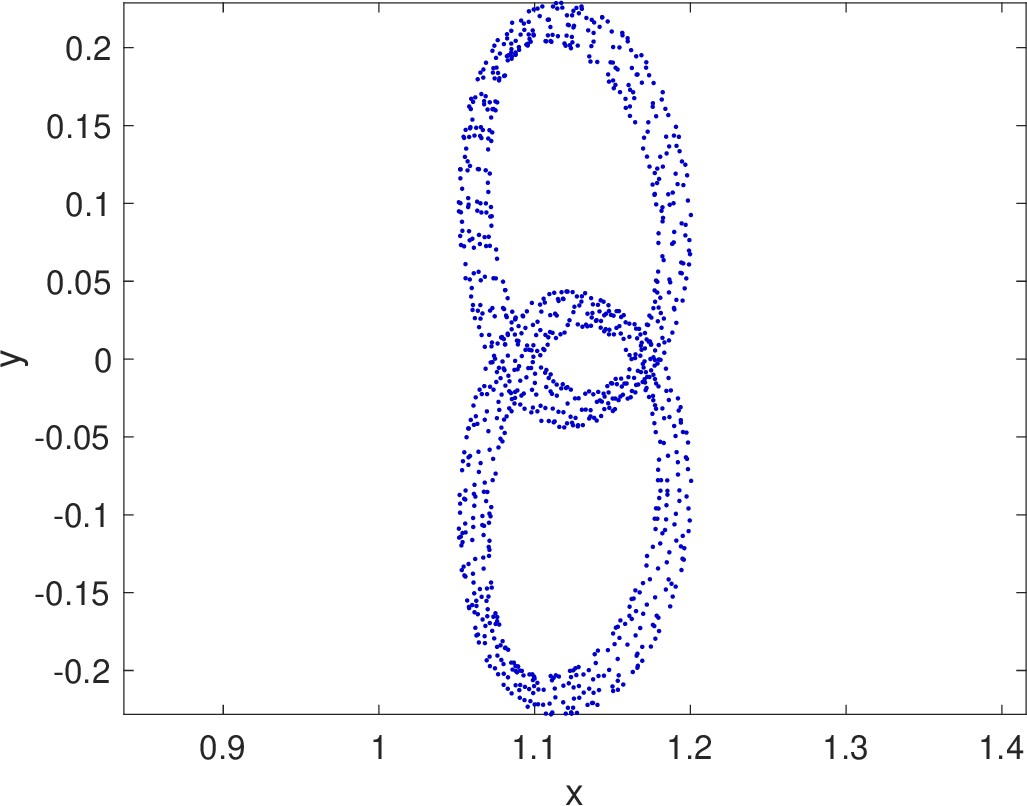}
}
\caption{
Quasiperiodic orbit tracked starting from base point $(1.13, -0.0902)$ at 
C = 3.14 in the constant rotating, non-pulsating frame.
}
\label{constantrotall}
\end{figure}

It is also interesting to examine the orbit in the constant rotating, non-pulsating frame
equivalent to the CRTBP's rotating frame.  The orbit and $\Sigma_{y=0}$ intersections
are plotted in this frame in
Figure \ref{constantrotall}, and some interesting characteristics may be observed.
The intersections in Figure \ref{constantrot} now show some structure that was
not apparent in the variable rotating, non-pulsating frame.  As expected, the
orbit also wanders further in the $y$ direction as the the frames are less aligned.

\clearpage

\section{Conclusions}

A method for performing isolating neighborhood computations in non-autonomous
systems was developed and applied to both a simple example and the ERTBP.  In
the case of the simple example, it was possible to compute more stringent
isolating block boundaries across all times, and in combination with a
bisection method, asymptotic target trajectories that would stay within the isolating
block for different initial times.  In the ERTBP, it was shown that it was possible to
compute isolating neighborhood boundaries that could be verified across a range
of energies.  By verifying that trajectories used in our orbit tracking
algorithm did not go outside of this range of energies, it was found that
quasiperiodic orbits around the L$_2$ libration point could be closely tracked
in the ERTBP.  Using this approach in the planar ERTBP, it was possible to
compute quasiperiodic orbits, or 2-tori, equivalent to the periodic planar Lyapunov orbit
in the CRTBP, and these results were found to be consistent with other results
in the literature that had been found for small eccentricities.  The method was
successfully extended to the spatial ERTBP, and orbits equivalent to the
Lissajous and quasi-halo orbits in the CRTBP were followed.  In these case, these orbits
had the characteristics of 3-tori, corresponding to the 2-tori in the CRTBP.
Different types of behavior were observed for the quasi-halo orbits in a small
region depending on the specific initial base point and initial time.
The developed algorithms have now been applied in the CRTBP and ERTBP, and this
approach lays the foundation for computations in other non-autonomous systems
such as the bicircular and ephemeris models.  This approach also allows for a relatively
seamless application to these other systems with perturbations that avoids many
of the difficulties inherent to other methods such as normal form or
parameterization approaches.



%

\section{Acknowledgements}

Part of the research presented in this paper has been carried out at the Jet
Propulsion Laboratory, California Institute of Technology, under a contract
with the National Aeronautics and Space Administration (80NM0018D0004).  
This work has been supported through funding by the Multimission Ground System
and Services Office (MGSS) in support of the development of the Advanced
Multi-Mission Operations System (AMMOS).

\bibliographystyle{AAS_publication}   
\bibliography{out}   

\clearpage
\noindent \textbf{APPENDIX}
\appendix

\subsection{A.  Non-Uniformly Rotating ERTBP Coordinate System}
\label{timeertbp}

An alternative formulation of the ERTBP that is useful for providing insight
into the problem is one where the $x$ axis rotates non-uniformly and stays
aligned with the two primaries.  Given this coordinate frame, the primary and
secondary will still be allowed to move along the $x$ axis over time.  Hiday
and Howell\cite{Hiday:1992b} define the equations of motion for this system
using time as the independent variable with their notation as
\begin{equation}
  \begin{aligned}
  \ddot{x} - 2 n \dot{y} &= \frac{\partial U}{ \partial x} + \dot{n}y \\
  \ddot{y} + 2 n \dot{x} &= \frac{\partial U}{ \partial y} - \dot{n}x \\
  \ddot{z} &= \frac{\partial U}{ \partial z}
  \end{aligned}
\end{equation}
where the pseudo-potential $U$ is given by
\begin{equation}
U = \frac{1}{2} n^2 (x^2 + y^2) + \frac{1-\mu}{r_1} + \frac{\mu}{r_2}
.
\end{equation}
Here, the variable angular velocity of the primaries is
\begin{equation}
n = \frac{\sqrt{1-e^2} }{(1 - e\cos E)^2}
\end{equation}
and
\begin{equation}
\dot{n} = -2e \frac{\sqrt{1-e^2} }{(1 - e\cos E)^4} \sin E
.
\end{equation}

\subsection{B.  Conversion Between Time and True Anomaly}

It is often convenient to convert between time and true anomaly, especially for
comparing integrated trajectories given using the current ERTBP formulation and
other formulations such as the one in Appendix A.  Standard conversions may be used
for this process.

\subsubsection{Conversion From Time to True Anomaly:}

If time is the independent variable, Kepler's equation
\begin{equation}
\label{keplerseqn}
  M = n (t - T) = E - e \sin E
\end{equation}
may be used to convert to $\nu$ where $T$ is the time of periapse crossing.
In this case the mean motion, $n$ = 1.
A Newton's method may be used to obtain $E$ from $M$.
Then, to obtain $\nu$,
\begin{equation}
  \nu = \mbox{acos}\left(\frac{e - \cos(E)}{e \cos(E) - 1}\right)
\end{equation}
while keeping in mind the appropriate quadrant checks.

\subsubsection{Conversion From True Anomaly to Time:}

In this case, the eccentric anomaly may be computed from
\begin{equation}
E = \mathrm{acos}\left(\frac{e+\cos\left(\nu \right)}{e\,\cos\left(\nu \right)+1}\right)
\end{equation}
while performing the appropriate quadrant checks.  Then the change in time may be found
using Equation \ref{keplerseqn}.

\end{document}